\documentclass{amsart}
\usepackage{amsmath}
\usepackage{amsxtra}
\usepackage{amscd}
\usepackage{amsthm}
\usepackage{amsfonts}
\usepackage{amssymb}
\usepackage{eucal}

\newtheorem{lem}{Lemma}
\newtheorem{prop}{Proposition}

\newtheorem{thm}{Theorem}

\newcommand{\thmref}[1]{Theorem~\ref{#1}}
\newcommand{\secref}[1]{Sect.~\ref{#1}}
\newcommand{\lemref}[1]{Lemma~\ref{#1}}

\newcommand{\propref}[1]{Proposition~\ref{#1}}

\newcommand{\nc}{\newcommand}
\nc{\ssec}{\subsection}
\nc{\sssec}{\subsubsection}
\nc{\on}{\operatorname}
\nc{\bi}{\bibitem}

\nc{\ZZ}{{\mathbb Z}}
\nc{\FF}{{\mathbf F}}
\nc{\CC}{{\mathbb C}}
\nc{\GG}{{\mathbb G}}
\nc{\PP}{{\mathbb P}}
\nc\one{{\mathbf 1}}
\nc{\DD}{{\mathbb D}}
\nc{\HH}{{\mathbf H}}

\nc{\Fq}{{\mathbb F}_q}
\nc{\Fqb}{\overline{\mathbb F}_q}
\nc{\Ql}{\ol{\mathbb Q}_\ell}

\nc\Gr{\on{Gr}}
\nc\Hom{\on{Hom}}
\nc\Conv{\on{Conv}}
\nc\Fl{\on{Fl}}
\nc\Iw{\on{I}}
\nc\Aut{\on{Aut}}
\nc\IC{\on{IC}}

\nc{\OO}{{\mathcal O}}
\nc{\T}{{\mathcal T}}
\nc{\F}{{\mathcal F}}
\nc{\E}{{\mathcal E}}
\nc{\D}{{\mathcal D}}
\renewcommand{\H}{{\mathcal H}}
\nc{\X}{{\mathcal X}}
\nc{\U}{{\mathcal U}}
\nc{\V}{{\mathcal V}}
\nc{\W}{{\mathcal W}}
\nc{\Y}{{\mathcal Y}}
\renewcommand{\L}{{\mathcal L}}
\nc{\G}{{\mathcal G}}
\nc{\C}{{\mathcal C}}
\renewcommand{\S}{{\mathcal S}}
\renewcommand\O{\widehat{\OO}}
\nc\K{\widehat{\mathcal K}}

\nc\starfl{\underset{\Fl}\star}
\nc\stargr{\underset{\Gr}\star}

\nc\ol{\overline}
\nc\wt{\widetilde}
\nc\tboxtimes{\wt{\boxtimes}}

\title[Central elements in the affine Hecke algebra via nearby cycles]
{Construction of central elements in the affine Hecke \\
algebra via nearby cycles}

\author{D.~Gaitsgory}

\address{Department of Mathematics, Harvard University \\
Cambridge MA 02138}

\email{gaitsgde@math.harvard.edu}

\begin{document}

\maketitle

\section*{Introduction}

\ssec{Overview}

Let $G$ be a connected reductive group over a finite field $\Fq$ and let $G(\K)$ be the corresponding group
over the local field $\K=\Fq((t))$. Let $G(\O)\subset G(\K)$ be a maximal compact subgroup 
of $G(\K)$ (here $\O=\Fq[[t]]$) and let $\HH_{sph}$ denote the Hecke algebra of $G(\K)$ with 
respect to $G(\O)$.

In other words, $\HH_{sph}$ as a vector space consists of compactly supported 
bi--$G(\O)$--invariant functions
$G(\K)\to\Ql$ and the product is defined by 
$$f^1\star f^2(g)=\underset{G(\K)}\int f^1(x)\cdot f^2(x^{-1}\cdot g) \; dx,$$
where $dx$ is a Haar measure on $G(\K)$ with $dx(G(\O))=1$.

The basic fact about $\HH_{sph}$ is that it is commutative. Moreover, one can show that when $G$
is split, $\HH_{sph}$
is isomorphic to the Grothendieck ring of the category of finite dimensional representations of the
Langlands dual group $\check G$.

\medskip

Now, let $\Iw\subset G(\O)$ be the Iwahori subgroup, and let us consider the corresponding Hecke 
algebra,
denoted $\HH_{\Iw}$. For example, when $G$ is simply--connected, $\HH_{\Iw}$ can be identified 
with the affine Hecke algebra attached to the root system of $G$.

Unlike $\HH_{sph}$, the algebra $\HH_{\Iw}$ is non--commutative, and in this paper we will 
be concerned
with its center, denoted $Z(\HH_{\Iw})$. The starting point is a theorem saying
that $Z(\HH_{\Iw})\simeq \HH_{sph}$. Moreover, the map in one direction can be described very 
explicitly.

Let $\pi$ be a linear map from $\HH_{\Iw}$ to the space of $G(\O)-\Iw$--invariant functions 
defined by 
$$\pi(f)(g)=\underset{\Iw\backslash G(\O)}\int f(x\cdot g) \; dx.$$
It is easy to see that $\pi$ maps $Z(\HH_{\Iw})$ to $\HH_{sph}$ and a theorem of 
J.~Bernstein (cf. \cite{Be}, Theorem 2.13 or \cite{Lu}, Proposition 8.6)
asserts that (at least when $G$ is split) this is an isomorphism. Our goal
in this paper is to describe in some sense explicitly the inverse map to $\pi$. 
This will be done by realizing $\HH_{sph}$ and $\HH_{\Iw}$ geometrically. 

\bigskip

First, there exists a group--scheme (resp., an group--indscheme) over $\Fq$ whose set of 
$\Fq$--points 
identifies with $G(\O)$ (resp., with $G(\K)$). We will abuse the notation and denote these 
objects again
by $G(\O)$ and $G(\K)$, respectively. In addition, there exists a subgroup $\Iw\subset G(\O)$ 
of finite codimension,
such that the quotient $G(\O)/\Iw$ is the flag variety $G/B$, where $B$ is a Borel subgroup of $G$.

One can form the quotients $\Gr=G(\K)/G(\O)$ and $\Fl=G(\K)/\Iw$, which will be indschemes over 
$\Fq$ and study the categories of perverse sheaves: $\on{P}_{G(\O)}(\Gr)$ 
(resp., $\on{P}_{\Iw}(\Fl)$)
will stand for the category of $G(\O)$--equivariant perverse sheaves on $\Gr$  
(resp., for $\Iw$--equivariant perverse sheaves on $\Fl$). 

The ``faisceaux--fonctions'' correspondence gives a map from the Grothendieck group of 
$\on{P}_{G(\O)}(\Gr)$
to $\HH_{sph}$ and from the Grothendieck group of $\on{P}_{\Iw}(\Fl)$ to $\HH_{\Iw}$. 
Moreover, one can introduce
convolution functors $\on{P}_{G(\O)}(\Gr)\stargr\on{P}_{G(\O)}(\Gr)\mapsto \on{P}_{G(\O)}(\Gr)$
and $\on{P}_{\Iw}(\Fl)\starfl\on{P}_{\Iw}(\Fl)\mapsto \on{D}^b_{\Iw}(\Fl)$ that will lift 
the $\star$
operations on $\HH_{sph}$ and $\HH_{\Iw}$, respectively (cf. \secref{grfl} for more details).

\medskip

Now, we can formulate our task more precisely: we would like to construct a functor
$Z: \on{P}_{G(\O)}(\Gr)\to \on{P}_{\Iw}(\Fl)$, such that on the level of Grothendieck groups
it induces the map $\pi^{-1}$.

It will turn out that this functor indeed exists and can be constructed using the operation
of taking nearby cycles of a perverse sheaf. Namely, we will construct a $1$--parametric family 
of schemes, which we will call $\Fl_X$, which degenerates the product $\Gr\times G/B$ to $\Fl$. 
Then for $\S\in \on{P}_{G(\O)}(\Gr)$, $Z(\S)$ will be the nearby cycles of the product
$\S\boxtimes \delta_{1_{G/B}}$.

Moreover, it will turn out that the functor $Z$ has some extremely favorable properties 
(cf. formulation
of \thmref{main}). In addition, since $Z$ is obtained by a nearby cycles construction, 
the perverse sheaves $Z(\S)$ will possess an extra structure: that of a nilpotent endomorphism, 
coming
from the monodromy. This phenomenon is invisible on the classical level (i.e. when one looks
at the corresponding Grothendieck groups and not at the categories), and supposedly it carries 
a deep representation--theoretic meaning (\cite{Bez}).

\medskip

\ssec{Conventions}  \label{conventions}

This paper uses in an extensive way the language of indschemes and of perverse sheaves on them. Although 
the objects we will operate with are straightforward extensions of the corresponding finite-dimensional ones,
not all of the definitions are present in the published literature, and for the reader's convenience
we will review them in the Appendix, \secref{App}.

\medskip

As was mentioned before, $G$ is a connected reductive group over the base field $\Fq$. By 
$\on{Rep}(G)$ we will denote the category of finite-dimensional $G$--representations. Throughout the paper, the 
notation 
$\F_G$ is reserved for principal $G$--bundles on various schemes and $\F^0_G$ we will denote {\it the} trivial 
$G$--bundle.

\smallskip

In several places in this paper we will use the concepts of a formal disc 
$\D$ and of a formal punctured disc $\D^*$.
They will appear in the following circumstances:

Let $S=\on{Spec}(\OO_S)$ be an affine scheme. An $S$--family of $G$--bundles on $\D$ 
(resp., on $\D^*$) is by definition a tensor functor from $\on{Rep}(G)$ 
to the tensor category of $S$--families of vector bundles on $\D$ (resp., $\D^*$), 
where the latter consists of finitely generated projective modules
over $\OO_S[[t]]$ (resp., $\OO_S((t))$). 

Let $\D_k=\on{Spec}(\Fq[[t]]/t^{k+1})$. It is easy to see that an $S$--family of $G$--bundles on
$\D$ is the same as a compatible system of $G$--bundles on $\D_k\times S$. If $\F_G$ is an $S$--family  
on $\D$, we will denote by $\F_G|_{\D^*}$ (resp., $\F_G|_{\D_k}$) the corresponding induced family on $\D^*$
(resp., on $\D_k$).

\smallskip

An $S$--family of maps $\D\to G$ (resp., $\D^*\to G$, $\D_k\to G$) is 
is by definition a ring homomorphism $\OO_G\to\OO_S[[t]]$ 
(resp., $\OO_G\to \OO_S((t))$, $\OO_G\to \OO_S[[t]]/t^{k+1}$). The functor that attaches to $S$ the set of all 
$S$--families of maps $\D\to G$ (resp., $\D^*\to G$, $\D_k\to G$) 
is representable by a group--scheme (resp., by an group--indscheme, algebraic group) 
that will be denoted $G(\O)$ (resp., $G(\K)$, $G(\O)_k$). We have: $G(\O)=\underset{\longleftarrow}{G(\O)_k}$.

\medskip

Let $\Y$ be a scheme, $H$ be an algebraic group and $\Y_1$ an $H$--torsor
over $\Y$. Let, in addition, $\Y_2$ be an $H$--scheme.
We will denote by $\Y_1\overset{H}\times \Y_2$ the associated fibration over $\Y$.
If $\T$ is a perverse sheaf on $\Y$ and $\S$ is an
$H$--equivariant perverse sheaf on $\Y_2$, we can form their
twisted external product $\T\tboxtimes\S$, which will be a perverse sheaf 
$\Y_1\overset{H}\times \Y_2$.

For a scheme $\Y$, $\Ql{}_\Y$ will denote the constant sheaf on $\Y$ and for $y\in\Y$, 
$\delta_y$ will denote the corresponding $\delta$--function sheaf.

\medskip

Finally, we should mention that although we work over the ground field $\Fq$, all the results
of this paper are valid over a ground field of characteristic zero.

\ssec{Acknowledgments}

This paper owes its existence to A.~Beilinson: the very idea of obtaining ``central''
objects of $\on{P}_{\Iw}(\Fl)$ as nearby cycles of objects of $\on{P}_{G(\O)}(\Gr)$ 
is an invention of his. 
\footnote{As Beilinson points out, he was in turn inspired by T.~Haines and R.~Kottwitz,
who proposed a similar idea in the framework of Shimura varieties, which has been
realized in a recent preprint by T.~Haines and B.~C.~Ngo.} 
The author wishes to thank R.~Bezrukavnikov
for stimulating discussions. Finally, I am grateful to Ya.~Varshavsky and the referee who pointed
out numerous mistakes  and whose comments helped me to 
improve the exposition. 

\section{Formulation of the results}

\ssec{Affine Grassmannian and affine flags} \label{grfl}

\sssec{}

Consider the functor that associates to a scheme $S$ the set of pairs
$(\F_G,\beta)$, where $\F_G$ is an $S$--family of $G$--bundles on $\D$
and $\beta$ is a trivialization of the corresponding family of $G$--bundles on $\D^*$, i.e.
$\beta:\F_G|_{\D^*}\to \F^0_G|_{\D^*}$. This functor is representable by an indscheme 
(cf. \secref{App}), which we will denote by $\Gr$,
called the affine Grassmannian of $G$. 

Here are the basic properties of $\Gr$. First, $\Gr$ has a distinguished point 
$1_{\Gr}\in\Gr$ that 
corresponds to the pair $(\F^0_G,\beta^0)$, where $\beta^0$ 
is the tautological trivialization of the trivial bundle.

Consider the group--scheme $G(\O)$ and the group--indscheme 
$G(\K)$ (cf. \secref{conventions}). It is obvious that for a scheme $S$, $\on{Hom}(S,G(\K))$ 
is the group of automorphisms
of the trivial $S$--family of $G$--bundles on $\D^*$. Hence, $G(\K)$ acts on $\Gr$ in a natural way, by changing the
data of $\beta$. 

It is known that the induced $G(\O)$--action on $\Gr$ is ``nice'' (cf. \secref{App}). 
This means that $\Gr$ can be represented 
as a union of finite-dimensional closed subschemes, each of which is $G(\O)$-stable. 
Therefore, we can introduce the 
category $\on{P}_{G(\O)}(\Gr)$ of $G(\O)$--equivariant perverse sheaves on $\Gr$
along with the corresponding derived category $\on{D}^b_{G(\O)}(\Gr)$.

\sssec{}

Now we will recall the convolution operation 
$\on{P}_{G(\O)}(\Gr)\times \on{P}(\Gr)\mapsto \on{D}^b(\Gr)$.

For a non-negative integer $k$, let $\G_k$ be a $G(\O)_k$--torsor over $\Gr$, defined as 
the indscheme that represents the functor that associates to a scheme $S$ a triple
$(\F_G,\beta,\gamma_k)$, where $(\F_G,\beta)$ are as above and $\gamma_k$ is a 
trivialization of $\F_G|_{\D_k}$. 

This defines a $G(\O)$--torsor $\G$ over $\Gr$, cf. \secref{App}. If we were to consider the 
total space of $\G$, it would be an indscheme (not of ind--finite type), 
isomorphic to $G(\K)$. Therefore, one can loosely speak of $\Gr$ as being the 
quotient $G(\K)/G(\O)$, which we will sometimes do in order to save notation.

The convolution diagram, denoted $\Conv_{\Gr}$, is the indscheme, associated to the $G(\O)$--torsor $\G$ over $\Gr$
and the $G(\O)$--scheme $\Gr$, i.e. $\Conv_{\Gr}=\G\overset{G(\O)}\times \Gr$, according to our conventions. We
again refer the reader to \secref{App} for the explanation why $\Conv_{\Gr}$ is a well-defined indscheme, as
well as for the proof of the following lemma:

\begin{lem}  \label{functconv}
The indscheme $\Conv_{\Gr}$ represents the functor that attaches to a scheme $S$ a quadruple
$(\F_G,\F^1_G,\wt{\beta},\beta^1)$, where $\F_G$ and $\F^1_G$ are $S$--families of $G$--bundles on $\D$,
$\wt{\beta}$ is an isomorphism $\F_G|_{\D^*}\to\F^1_G|_{\D^*}$ between the induced families of $G$--bundles
on $\D^*$ and $\beta^1$ is an isomorphism $\F^1_G|_{\D^*}\to\F^0_G|_{\D^*}$.
\end{lem}

There are two natural projections $p,p':\Conv_{\Gr}\to\Gr$. In the above functorial language,
$p^1$ sends a quadruple $(\F_G,\F^1_G,\wt{\beta},\beta^1)$ to $(\F^1_G,\beta^1)$ and 
$p(\F_G,\F^1_G,\wt{\beta},\beta^1)=(\F_G,\beta^1\circ\wt{\beta})$. Naively, one should picture the above
projections as follows: if we identify $\Conv_{\Gr}$ with $G(\K)\overset{G(\O)}\times G(\K)/G(\O)$, then
$p^1$ is the projection on the first factor, i.e. $p^1(g_1\times g)=g_1$ and $p(g_1\times g)=g_1\cdot g$.

\smallskip

Thus, $p^1$ realizes $\Conv_{\Gr}$ as a fibration over $\Gr$, with the typical 
fiber isomorphic again to $\Gr$. We are going to use the twisted external product
construction, introduced in \secref{conventions} and extended for ind-schemes
in \secref{App}:

Starting with an object $\T\in \on{P}(\Gr)$ and an object 
$\S\in \on{P}_{G(\O)}(\Gr)$, we can form their twisted external product $\T\tboxtimes\S$, 
which will be an object of $\on{P}(\Conv_{\Gr})$.

Finally we set $\T\stargr\S:=p_{!}(\T\tboxtimes\S)\in \on{D}^b(\Gr)$. 
It is easy to see that if $\T$ is also an object of 
$\on{P}_{G(\O)}(\Gr)$, then $\T\stargr\S$ will belong to $\on{D}^b_{G(\O)}(\Gr)$.
Thus, $\stargr$ induces a bi--functor from $\on{P}_{G(\O)}(\Gr)$ to 
$\on{D}^b_{G(\O)}(\Gr)$. 
On the level of Grothendieck groups, $\stargr$ descends, of course, to the usual 
convolution product on the spherical Hecke algebra $\HH_{sph}$.

\smallskip

\noindent{\it Remark.} It 
follows from Lusztig's work \cite{Lu} that for $\S_1,\S_2\in \on{P}_{G(\O)}(\Gr)$, 
the convolution $\S_1\stargr\S_2$ is again a 
perverse sheaf, i.e. the $\stargr$--operation makes
$\on{P}_{G(\O)}(\Gr)$ into a monoidal category. 
Moreover, the fact that the spherical Hecke algebra
$\HH_{sph}$ is commutative can be lifted to the categorical level: 
one can endow $\on{P}_{G(\O)}(\Gr)$
with a commutativity constraint, 
i.e. $\on{P}_{G(\O)}(\Gr)$ has a structure of a tensor category.

As a by--product of the results of this paper, we will construct the commutativity constraint
$\S_1\stargr\S_2\to \S_2\stargr\S_1$ and, 
in addition, we will prove a strengthened version of Lusztig's
theorem: we will show that for $\S\in \on{P}_{G(\O)}(\Gr)$,
$\S\tboxtimes\T$ is a perverse sheaf for any $\T\in \on{P}(\Gr)$. 

\sssec{}

Let us fix once and for all a Borel subgroup $B\subset G$, it corresponds to a 
distinguished point $1_{G/B}$
inside the flag variety $G/B$. The Iwahori group $\Iw\subset G(\O)$ is by definition the
preimage of $B$ under the natural projection $G(\O)\to G(\O)_0=G$. We will denote by 
$\Iw_k$ the image of $\Iw$ under the projection $G(\O)\to G(\O)_k$. We have: 
$\Iw=\underset{\leftarrow}{\Iw_k}$.

The affine flag variety $\Fl$ is the indscheme associated to the $G$--bundle $\G_0$ over $\Gr$ and a $G$--variety $G/B$,
i.e. $\Fl=\G_0\overset{G}\times G/B$. Thus, loosely speaking, $\Fl=G(\K)/\Iw$. We will denote by $\pi$
the natural projection $\Fl\to\Gr$.

\smallskip

Functorially, for a scheme $S$, the set $\on{Hom}(S,\Fl)$ consists of triples
$(\F_G,\beta,\epsilon)$, where $\F_G$ and $\beta$
are as in the definition of $\Gr$ and $\epsilon$ is a
reduction of $\F_G|_{\D_0}$ to $B$.

Let $1_{\Fl}\in\Fl$ be the distinguished point that corresponds to the triple $(\F^0_G,\beta^0,\epsilon^0)$,
where $(\F^0_G,\beta^0)=1_{\Gr}$ and $\epsilon^0$ corresponds to the chosen Borel subgroup $B\subset G$.

The $G(\K)$-action on $\Gr$ lifts in a natural way to an action on $\Fl$ and the 
induced actions of $G(\O)$ and hence of $\Iw$ are ``nice'' in the sense of \secref{App}. Therefore, one may
consider the categories $\on{P}(\Fl)$, $\on{D}^b(\Fl)$, $\on{P}_{\Iw}(\Fl)$ and $\on{D}_{\Iw}(\Fl)$. 

As in the case of $\Gr$, one defines the ind--scheme $\Conv_{\Fl}$, which classifies the data of 6-tuples
$(\F_G,\F^1_G,\wt{\beta},\beta^1,\epsilon,\epsilon^1)$, where $(\F_G,\F^1_G,\wt{\beta},\beta^1)$ are as in
the definition of $\Conv_{\Gr}$ and $\epsilon$ (resp., $\epsilon^1$) is a reduction of 
$\F_G|_{\D_0}$ (resp., of $\F^1_G|_{\D_0}$) to $B$.

Let $p^1$ and and $p$ denote the two projections from 
$\Conv_{\Fl}$ to $\Fl$. As in the previous case,
we obtain a functor 
$\T,\S\mapsto \T\tboxtimes\S$ from $\on{P}(\Fl)\times \on{P}_{\Iw}(\Fl)$
to $\on{D}^b(\Conv_{\Fl})$ and we set $\T\starfl\S:=p_!(\T\tboxtimes\S)$.

However, Lusztig's theorem does not extend to the case of affine 
flags: the convolution functor 
$\starfl$ does not preserve perversity, i.e. 
it maps $\on{P}_{\Iw}(\Fl)\times \on{P}_{\Iw}(\Fl)$
to $\on{D}_{\Iw}(\Fl)$. In addition, there certainly is no isomorphism
$\S_1\starfl\S_2\to \S_2\starfl\S_1$ since the corresponding equality is not true even on the Grothendieck group level 
(the Iwahori Hecke algebra $\HH_{\Iw}$ is not commutative).

\ssec{The functor $Z$}

\sssec{}

The main result of this paper is the following theorem:

\begin{thm}  \label{main}
There exists a functor: $Z:\on{P}_{G(\O)}(\Gr)\to \on{P}_{\Iw}(\Fl)$
possessing the following properties:

\smallskip

\noindent{\rm (a)} For $\S\in \on{P}_{G(\O)}(\Gr)$ and an
arbitrary perverse sheaf $\T$ on $\Fl$, the convolution
$\T\starfl Z(\S)$ is a perverse sheaf.

\smallskip

\noindent{\rm (b)} For $\S\in \on{P}_{G(\O)}(\Gr)$ and
$\T\in \on{P}_{\Iw}(\Fl)$ there is a canonical isomorphism
$Z(\S)\starfl \T\simeq \T\starfl Z(\S)$.

\smallskip

\noindent{\rm (c)} We have $Z(\delta_{1_{\Gr}})=\delta_{1_{\Fl}}$ and
for $\S^1,\S^2\in \on{P}_{G(\O)}(\Gr)$
there is a canonical isomorphism
$Z(\S^1)\starfl Z(\S^2)\simeq Z(\S^1\stargr\S^2)$.

\smallskip

\noindent{\rm (d)}
For $\S\in \on{P}_{G(\O)}(\Gr)$, we have $\pi_{!}(Z(\S))\simeq \S$.

\end{thm}

\medskip

Property (d) above insures that the composition 
$$K(\on{P}_{G(\O)}(\Gr))\overset{Z}\to K(\on{P}_{\Iw}(\Fl))\to \HH_{\Iw}$$ equals
$K(\on{P}_{G(\O)}(\Gr))\to \HH_{sph}\simeq Z(\HH_{\Iw})$, i.e.
\thmref{main} fulfills our promise to construct geometrically the inverse of the map
$Z(\HH_{\Iw})\to \HH_{sph}$.

\sssec{}

An additional basic structure of the functor $Z$ is described by the following theorem:

\begin{thm} \label{monodromy}

The functor $Z$ carries a nilpotent (monodromy) endomorphism $M$ 
$$M_\S: Z(\S)\to Z(\S)(-1), \text{ for } \S\in \on{P}_{G(\O)}(\Gr),$$
which is compatible with the isomorphisms 
of \thmref{main}(c): For $\S_1,\S_2\in \on{P}_{G(\O)}(\Gr)$ the
square
$$
\CD
Z(\S^1)\starfl Z(\S^2) 
@>{M_{\S^1}\starfl\on{id}_{\S^2}+\on{id}_{\S^1}\starfl M_{\S^2}}>> 
Z(\S^1)\starfl Z(\S^2)(-1)  \\
@VVV            @VVV   \\
Z(\S^1\stargr\S^2)   @>{M_{\S^1\stargr\S^2}}>>   Z(\S^1\stargr\S^2)(-1)
\endCD
$$
commutes. 

\end{thm}

\sssec{}

Let explain the 
first non-trivial example of how \thmref{main} works for $G=GL(2)$.
Consider the following closed $G(\O)$--stable subscheme $Y_0$ of $\Gr$:

By definition, for $G=GL(2)$, $\Gr$ classifies lattices in $\K\oplus \K$ (i.e. 
$\O$-submodules $\L\subset\K\oplus \K$ of rank $2$) and $Y_0$ corresponds
to those $\L$ which are contained in $\L^0:=\O\oplus\O$ with $\dim(\L^0/\L)=1$.

By construction, $Y_0$ is isomorphic to the projective line ${\mathbb P}^1$.
We take $\S\in \on{P}_{G(\O)}(\Gr)$ to be $\Ql{}_{Y_0}[1]$. In this case, the construction which will be discussed in the 
next section reduces to the usual Picard-Lefschetz situation
and the perverse sheaf $Z(\S)\in \on{P}_{\Iw}(\Fl)$ can be described very explicitly.

First, $Z(\S)$ will be supported on $\pi^{-1}(Y_0)\subset\Fl$. Now let $Y_1$ and $Y_2$
be the following two subschemes of $\pi^{-1}(Y_0)$: 

By definition, $\Fl$ classifies lattices $\L$ as above plus a choice of a line 
$\ell\subset\L/t\cdot\L$, where $t$ is the uniformizer of $\O$. The chosen Borel subgroup 
$B\subset GL(2)$ fixes a line $\ell^0\subset \L^0/t\cdot \L^0\simeq\Fq\oplus\Fq$
and $Y_1$ corresponds to the lattice $L^1=\on{ker}(\L_0\to \Fq\oplus\Fq\to \Fq\oplus\Fq/\ell^0)$
and an arbitrary $\ell$. On the contrary, $Y_2$ corresponds to an 
arbitrary $\L\in Y_0$, but $\ell$ must be the kernel of the map 
$\L/t\cdot\L\to  \L^0/t\cdot \L^0$. Both $Y_1$ and $Y_2$
are isomorphic to ${\mathbb P}^1$ and their intersection in $\Fl$ 
is a point-scheme, which we will denote by $Y_3$.

We claim that $Z(\S)$ has the following form: it has a three-step filtration 
$$0=F_0\subset F_1\subset F_2\subset F_3=Z(\S),$$ such that $F_1\simeq \delta_{Y_3}$,
$F_2/F_1\simeq \Ql{}_{Y_1}[1]\oplus \Ql{}_{Y_2}[1]$ and $F_3/F_2\simeq \delta_{Y_3}(-1)$. Moreover,
the monodromy map $M_{\S}$ acts as follows:
$$Z(\S)\twoheadrightarrow F_3/F_2\simeq \delta_{Y_3}(-1)\simeq F_1(-1)\hookrightarrow Z(\S)(-1).$$

\section{Construction-I}

\ssec{Global version of the affine Grassmannian}

\sssec{}

To carry out our constructions we will need to choose a curve $X$, which is smooth,
but not necessarily complete, and an $\Fq$--point point $x\in X$. In what follows, we will choose once
and for all an identification of the completed local ring $\O_x$ with $\O$.

The starting point is the following result of \cite{BL}:

\begin{lem} \label{bl}

The indscheme $\Gr$ represents the following functor: $\Hom(S,\Gr)$ is the set of pairs
$(\F_G,\beta)$, where $\F_G$ is a $G$--bundle on $X\times S$ and $\beta$ is its trivialization
over $(X\setminus x)\times S$.

\end{lem}

\smallskip

\noindent{\it Remark.}  Of course, if one has a pair $(\F_G,\beta)$ as in the proposition, one can restrict 
it to the formal disc around $x$ and thus obtain a point of $\Gr$ in the original definition. The meaning 
of \lemref{bl} is that this restriction is a bijection between the data on $X$ and that on the formal disc.

\sssec{}

Let $\Aut$ denote the pro--algebraic group of automorphisms of $\O$, i.e. for an affine 
scheme $S$, $\on{Hom}(S,\Aut)$ consists of all $\OO_S$--linear continuous automorphisms of 
$\OO_S[[t]]$. By definition, $\Aut$ is the projective limit of the groups $\Aut_k$,
where $\Aut_k$ is the group of automorphisms of $\Fq[[t]]/t^{k+1}$.
This group acts in a 
canonical way on $G(\O)$, $G(\K)$ and on $\Gr$; moreover, its action on $\Gr$ is ``nice'' 
in the sense of \secref{App}.

In addition, there is a canonical $\Aut$--torsor $\X$ over $X$: for an affine scheme $S$, 
an $S$--point of $\X$ is a pair $y:S\to X$ and a continuous $\OO_S$--linear 
isomorphism between $\OO_S[[t]]$ and the completion of
$\OO_{X\times S}$ along the graph $\Gamma_y$ of the map $y$.

We define the global version of $\Gr$ over $X$, denoted $\Gr_X$, 
as an indscheme associated to the $\Aut$--torsor $\X$ over $X$ and the 
$\Aut$--scheme $\Gr$, i.e. $\Gr_X:=\X\overset{\Aut}\times \Gr$,
cf. \secref{App}.

By invoking again the theorem of \cite{BL}, we obtain the following:

\begin{lem} \label{morebl}
The indscheme $\Gr_X$ represents the following functor: $\Hom(S,\Gr_X)$ is the set of triples
$(y,\F_G,\beta)$, where $y$ is an $S$--point of $X$, $\F_G$ is a $G$--bundle on $X\times S$ and $\beta$
is a trivialization $\F_G|_{X\times S\setminus\Gamma_y}\to\F^0_G|_{X\times S\setminus\Gamma_y}$,
where $\Gamma_y\subset X\times S$ is the graph of $y:S\to X$.
\end{lem}

We will denote by $\Gr_{X\setminus x}$ (resp., $\Gr_x$) the preimage of $X\setminus x$ (resp., of $x\in X$)
under the natural projection $\Gr_X\to X$.

\sssec{}  \label{globalsheaf}

An important observation is that to an object $\S\in \on{P}_{G(\O)}(\Gr)$, one can attach in a 
canonical way a perverse sheaf on $\Gr_X$. First, we have the following assertion:

\begin{prop} \label{GOeq}
Every $G(\O)$--equivariant perverse sheaf on $\Gr$ is automatically equivariant with respect to
$\Aut$.
\end{prop}

\begin{proof}

To prove the proposition, it is enough to show that over an algebraic closure $\Fqb$ of $\Fq$, the category
$\on{P}_{G(\O)}(\Gr)$ is semi--simple and every irreducible object in it is $\Aut$--equivariant. Indeed, this
would imply that every object of $\on{P}_{G(\O)}(\Gr)$ has a form $\underset{i}\oplus \,\S_i\otimes V_i$, where $\S_i$
is absolutely irreducible (and hence $\Aut$--equivariant) and $V_i$ is a ``perverse sheaf over $\on{Spec}(\Fq)$'',
i.e. a vector space acted on by the Frobenius.

The two facts mentioned above are well-known and we include the proof for completeness.
First of all, every irreducible $G(\O)$--equivariant perverse 
sheaf on $\Gr$ is an intersection 
cohomology sheaf on a closure of a $G(\O)$--orbit and each such orbit
is $\Aut$-stable. 

To prove the semi-simplicity assertion, we must show that if 
$Y'$ and $Y$ are two $G(\O)$--orbits, then  
$\on{Ext}^1(\on{IC}_{\overline{Y}},\on{IC}_{\overline{Y'}})=0$. 

First, let us assume that $Y=Y'$. The natural map 
$$\on{Ext}^1(\on{IC}_{\overline{Y}},\on{IC}_{\overline{Y'}})
\to \on{Ext}^1(\Ql{}_Y,\Ql{}_Y)$$
is an injection. However, the right-hand side is nothing but $H^1(Y,\Ql)$, and it vanishes,
since any $G(\O)$--orbit $Y$ is a isomorphic to a principal bundle over a (partial) flag variety
of the group $G$ with a unipotent structure group.

Thus, let $Y\neq Y'$ and without restricting the generality, we can assume that $Y'\subset\overline{Y}$. 
(It is easy to see that when neither $Y'\subset\overline{Y}$ nor $Y\subset\overline{Y'}$, the above
$\on{Ext}^1$ is automaticaly zero). It is enough to show that the $*$-restriction
of $\on{IC}_{\overline{Y}}$ to $Y'$ lives in the perverse cohomological degrees $\leq -2$.
It is known due to \cite{KL}, that the stalks of $\on{IC}_{\overline{Y}}$
have the parity vanishing property. \footnote{Following the referee's suggestion,
since the last assertion was not stated explicitly for the affine case in \cite{KL},
we will give a self-conated proof at the end of \secref{App}.}
Since the dimensions of $Y$ and $Y'$
have the same parity, $\on{IC}_{\overline{Y}}|_{Y'}$ has perverse cohomologies only in the even
degrees. In particular, its $-1$-st perverse cohomology sheaf is $0$, which is what
we had to prove. \footnote{The above proof that 
$\on{Ext}^1(\on{IC}_{\overline{Y}},\on{IC}_{\overline{Y'}})=0$ over
an algebraically closed field relies on the Kazhdan-Lusztig
parity vanishing assertion. However, in the recent paper \cite{FGV} it was shown that the fact
that $\on{Ext}^1(\on{IC}_{\overline{Y}},\on{IC}_{\overline{Y'}})=0$ for $Y=Y'$ formally implies 
the vanishing in the general case.}

\end{proof}

Thus, starting with $\S\in \on{P}_{G(\O)}(\Gr)$, or, 
more generally, with $\S\in \on{P}_{\Aut}(\Gr)$,
we can attach to it a perverse sheaf $\S_X\in\on{P}(\Gr_X)$,
by taking the twisted external product with the constant perverse sheaf $\Ql[1]$ on $X$.
We will denote the restriction of $\S_X$ to $\Gr_{X\setminus x}$ by $\S_{X\setminus x}$.

\sssec{}

The following will be useful in the sequel:

Starting from the $\Aut$--torsor $\X$ over $X$ and the group--scheme $G(\O)$, we can form a group--scheme
$G(\O)_X$ over $X$ by setting $G(\O)_X:=\X\overset{\Aut}\times G(\O)$. By construction, $G(\O)_X$
acts on $\Gr_X$ and for $\S\in \on{P}_{G(\O)}(\Gr)$ the perverse sheaf $\S_X$ 
are equivariant with respect to this action.

\begin{prop} \label{expequiv} 
Let $S$ be a scheme, let $(y,\F_G,\beta)$ be an $S$--point of $\Gr_X$ and let $\phi:X\times S\to G$ be a map.
Consider another $S$--point of $\Gr_X$ equal to $(y,\F_G,\phi\circ\beta)$. Then, for $\S\in \on{P}_{G(\O)}(\Gr)$
the pull-backs of $\S_X$ under these two maps from $\Gr_X$ to $S$ are canonically isomorphic
\end{prop}

\begin{proof}

By taking Taylor expansions of the map $\phi$, we obtain a map $\widehat{\phi}:X\times S\to G(\O)_X$.
Consider the map $S\to G(\O)_X\underset{X}\times\Gr_X$ obtained from the first map $S\to \Gr_X$
and 
$$S\overset{y\times\on{id}}\longrightarrow X\times S\overset{\widehat{\phi}}\to G(\O)_X.$$

Our two $S$--points of $\Gr_X$ are the compositions of the above map followed by the projection
$G(\O)_X\times \Gr_X\to \Gr_X$ in the first case, and by the action map $G(\O)_X\times \Gr_X\to \Gr_X$
in the second one.

Hence, the proposition follows from the $G(\O)_X$--equivariance of $\S_X$.

\end{proof}

\ssec{The nearby cycles construction}

\sssec{}

First, we will construct an indscheme $\Fl_X$ over $X$. We define $\Hom(S,\Fl_X)$ to be the set of
quadruples $(y,\F_G,\beta,\epsilon)$, where $(y,\F_G,\beta)$ are as in \lemref{morebl} and $\epsilon$
is a data of a reduction of $\F_G|_{x\times S}$ to $B$.

Obviously, $\Fl_X$ is a fibration over $\Gr_X$ with the typical fiber $G/B$. We will denote the projection
$\Fl_X\to\Gr_X$ by $\pi_X$. Let $\Fl_{X\setminus x}$ and $\Fl_x$ be the corresponding subschemes of $\Fl_X$.

\begin{prop} \label{descrFLX}

We have canonical isomorphisms $\Fl_{X\setminus x}\simeq \Gr_{X\setminus x}\times G/B$ and
$\Fl_x\simeq \Fl$.

\end{prop}

\begin{proof}

Let $(y,\F_G,\beta,\epsilon)$ be an $S$--point of $\Fl_X$ with $\Gamma_y\cap (x\times S)=\emptyset$.
Then the data of $\beta$ trivializes $\F_G|_{x\times S}$. Therefore, $\epsilon$ is
a reduction to $B$ of the trivial $G$--bundle on $S$, i.e. a map $S\to G/B$.

This defines a map $\Fl_{X\setminus x}\simeq \Gr_{X\setminus x}\times G/B$ and it is
straightforward to see that it is an isomorphism. 

The fact that $\Fl_x\simeq \Fl$ follows immediately from \lemref{bl}.

\end{proof}

\sssec{}  \label{nearby}

Let us recall the general formalism of the nearby cycles functor. Let $\Y$ be a scheme mapping 
to $X$
and let $\Y_{X\setminus x}$ and $\Y_x$ be its corresponding subschemes.

We have a functor 
$$\Psi_\Y:\on{D}^b(\Y_{X\setminus x})\to \on{D}^b(\Y_x),$$
whose basic property is that it maps $\on{P}(\Y_{X\setminus x})$ to $\on{P}(\Y_x)$, according to \cite{Bei}.
\footnote{{\it A priori}, for $\S\in \on{D}^b(\Y_{X\setminus x})$, $\Psi_\Y(\S)$
is defined only over $\Fqb$. To endow it with an $\Fq$--structure one needs to choose
a splitting from the Galois group $\on{Gal}(\Fqb/\Fq)$ to the Galois group 
of the field of fractions of the henselization of the local ring ${\mathcal O}_x$. From now on we choose
such a splitting.}

\sssec{} \label{intrconstr}

We apply $\Psi$ for $\Y=\Fl_X$. Using \secref{globalsheaf} and \lemref{descrFLX},
we can construct a functor from $\on{P}_{\Aut}(\Gr)$ to $\on{P}(\Fl_{X\setminus x})$:
$$\S\mapsto \S_{X\setminus x}\boxtimes \delta_{1_{G/B}}.$$

We set $Z(\S)=\Psi_{\Fl_X}(\S_{X\setminus x}\boxtimes \delta_{1_{G/B}})\in \on{P}(\Fl)$.

\medskip

It is straightforward to see that $Z(\delta_{1_{\Gr}})\simeq \delta_{1_{\Fl}}$. 
Indeed, we have a canonical
section $1_{\Fl_X}:X\to\Fl_X$ that sends $y$ the quadruple $(y,\F^0_G,\beta^0,\epsilon^0)$ and
$1_{\Fl_{X\setminus x}}=1_{\Gr_{X\setminus x}}\times 1_{G/B}$, $1_{\Fl_x}=1_{\Fl}$. 

\ssec{$\Iw$-equivariance}

\sssec{}

By construction, $\S\to Z(\S)$ is a functor between $\on{P}_{\Aut}(\Gr)$ and
$\on{P}(\Fl)$. The next proposition asserts that it defines a functor 
$\on{P}_{G(\O)}(\Gr)\to \on{P}_{\Iw}(\Fl)$, as required in \thmref{main}.

\begin{prop}  \label{Iwequiv}
For $\S\in \on{P}_{G(\O)}(\Gr)$, the perverse sheaf $Z(\S)$ on $\Fl$ is $\Iw$--equivariant.
\end{prop}

\begin{proof}

Let $Y$ be the support of $Z(\S)$ in $\Fl$. Choose an integer $k$ so that 
$\Iw$ acts on $Y$ through 
the quotient $\Iw\to\Iw_k$. Thus, we have to show that $Z(\S)$ is $\Iw_k$--equivariant.

Let $\Iw_{glob}$ be the sheaf of groups (on the category of all schemes
with Zarisky topology), which attaches to a scheme $S$ the
group of maps from the localization of $X\times S$ around $x\times S$ to $G$,
with the condition that $x\times S$ maps to $B$. By taking Taylor expansions at $x$, we obtain
a map of sheaves $\Iw_{glob}\to \underline{\Iw}$ 
(the underline means ``the sheaf represented by''). 

\begin{lem} \label{surj}
The composition
$\Iw_{glob}\to \underline{\Iw}\to\underline{\Iw_k}$
is a surjection of sheaves of groups.
\end{lem}

The proof will be given below. From this lemma we infer that it suffices to prove the following: 

\medskip

Let $S$ be a smooth scheme and let $U$ be an open subset in $X\times S$ containing $x\times S$. Let
$\phi:U\to G$ be a section of $\Iw_{glob}$. We have two maps from 
$S\times \Fl\to \Fl$:
one is the projection on the second factor and the other one is obtained 
by composing $S\to \Iw_{glob}\to \underline{\Iw}$ and the action
of $\Iw$ on $\Fl$. We must show that the pull-backs of $Z(\S)$ under these two 
maps from $\Fl$ to $S\times \Fl$ are isomorphic.

\medskip

From $\phi$ we obtain two maps 
$\phi^1,\phi^2:U\underset X\times \Fl_X\to \Fl_X$: the map $\phi^1$ is just the projection
on the second factor and the map $\phi^2$ is described as follows:

Let $S'$ be another scheme and let us consider an $S'$-point of $U\underset X\times \Fl_X$.
This amounts to having an $S'$-point $(y,\F_G,\beta,\epsilon)$ of $\Fl_X$
and a map $g:S'\to S$ such that the map $S'\overset{y\times g}\longrightarrow X\times S$ 
factors through $U\subset X\times S$. 

Let $U'$ denote the preimage of $U$ under the map 
$X\times S'\overset{\on{id}\times g}\longrightarrow X\times S$. Using $\phi$, we obtain a map $\phi':U'\to G$.
We need to produce another $S'$-point $(y,\F^\phi_G,\beta^\phi,\epsilon^\phi)$ of $\Fl_X$. 

By definition, $\F^\phi_G$ is the trivial bundle on $X\times S'\setminus \Gamma_y$ and is identified
with $\F_G$ over $U'$. We have: $U'\cup (X\times S\setminus \Gamma_y)=X\times S'$. Therefore, to define
$\F^\phi_G$ we need to define a gluing isomorphism $\F^0_G|_{U'\cap (X\times S'\setminus \Gamma_y)}\simeq
\F_G|_{U'\cap (X\times S'\setminus \Gamma_y)}$. The latter is obtained by composing $\beta$ and $\phi'$. 
The data of $\epsilon^\phi$ is by definition induced by $\epsilon$, via the identification
$\F^\phi_G|_{U'}\simeq \F_G|_{U'}$, since $U'$ contains $x\times S'$.
The data of $\beta^\phi$ follows from the construction. 

Both $\phi^1$ and $\phi^2$ are maps of schemes over $X$. Note that their values over $x$
are the two maps $S\times \Fl\to \Fl$ described above.
In addition, $\phi$ defines a map $\widehat{\phi}:U\to G(\O)_X$, and over $X\setminus x$
the maps $\phi^1$ and $\phi^2$ factor by means of $\widehat{\phi}$ through the action of $G(\O)_{X\setminus x}$ on
$\Fl_{X\setminus x}\simeq \Gr_{X\setminus x}\times G/B$, as in \propref{expequiv}.

\bigskip

Both $\phi^1$ and $\phi^2$ are smooth
and it is well-known that the functor of nearby cycles commutes with pull-backs
under smooth morphisms.
Therefore, it suffices to show that
$$\Psi_{U\underset X\times \Fl_X}(\phi^1{}^*(\S_{X\setminus x}\boxtimes \delta_{1_{G/B}}))\simeq
\Psi_{U\underset X\times \Fl_X}(\phi^2{}^*(\S_{X\setminus x}\boxtimes \delta_{1_{G/B}})).$$
However, as in \propref{expequiv} we obtain that  
$\phi^1{}^*(\S_{X\setminus x}\boxtimes \delta_{1_{G/B}})\simeq
\phi^2{}^*(\S_{X\setminus x}\boxtimes \delta_{1_{G/B}})$.

\end{proof}

\sssec{}

Finally, let us prove \lemref{surj}.

\begin{proof}

We will prove that if $S$ is affine,
the map $\Iw_{glob}(S)\to \Hom(S,\Iw_k)$ is a surjection.

The proof goes by induction. The assertion is obvious for $k=0$, since $I_0=B$. Therefore, it
suffices to show that any map $g_k:\D_k\times S\to G$ whose restriction to 
$\D_{k-1}\times S$ is trivial, can be lifted to a map from the localization
$(X\times S)_{x\times S}$ to $G$.

However, the group of maps $\{g_k\}$ as above is abelian and is isomorphic
to the group of functions on $S$ with values in the Lie algebra of $G$. 
In addition, we know that $\on{Lie}(G)$ is spanned
by $1$-dimensional subspaces that correspond to subgroups of $G$ isomorphic to
either $\GG_m$ or $\GG_a$. Hence, it suffices to analyze separately these two cases.

For $\GG_a$ the assertion is obvious: we can assume that $X$ is affine and then any function on $\D_k\times S$ 
can be extended to a function on $X\times S$. Similarly for $\GG_m$:

Let $t$ be a coordinate on $X$ around $x$.
Our map $g_k$ is a $k$-jet of a regular function on $X\times S$ of the form 
$1+t^k\cdot f$, where $f$ is a function on $S$. This $1+t^k\cdot f$
is the required map $(X\times S)_{x\times S}\to \GG_m$.

\end{proof}

Thus, we have constructed the functor $Z$. However, in order to prove \thmref{main}, we will
need to interpret the convolutions 
$Z(\S)\starfl\T$ and $\T\starfl Z(\S)$ too, in terms of nearby cycles.

\section{Construction-II}

\ssec{The indscheme $\Fl'_X$.}

\sssec{}

Our main tool will be the indscheme $\Fl'_X$ over $X$ defined as follows:

\smallskip

\noindent For a scheme $S$, 
$\Hom(S,\Fl'_X)$ is the set of quadruples $(y,\F_G,\beta',\epsilon)$, where 
$(y,\F_G,\epsilon)$ are as in the definition 
of $\Fl_X$, but $\beta'$ is now a trivialization of $\F_G$ 
off the divisor $\Gamma_y\cup (x\times S)$.

\medskip

Analogously, we introduce the indscheme $\Gr'_X$: $\Hom(S,\Gr'_X)$ is the set of triples 
$(y,\F_G,\beta')$, where $y$, $\F_G$ and $\beta'$ are as above. Of course, $\Fl'_X$ is a fibration over $\Gr'_X$
with the typical fiber $G/B$ and we will denote by $\pi'_X$ the corresponding projection.  Note that our
$\Gr'_X$ is a particular case of the Beilinson-Drinfeld Grassmannian, which was studied in \cite{MV} and \cite{BD}.

The fact that the above functors are indeed representable by indschemes can be proved by
a straightforward generalization of the argument that shows that $\Gr$ is representable 
by an indscheme (cf. \secref{App}). Here is a rough outline of the proof:

\begin{proof}(sketch)

For two integers $m,n\in {\mathbb N}$ we introduce a (relative) Hilbert scheme $\on{Hilb}^m_n$, where 
for a scheme $S$, $\on{Hom}(S,\on{Hilb}^m_n)$ consists of a map $y:S\to X$ and a coherent subsheaf
${\mathcal J}$ of $$\OO^{\oplus n}(m\cdot(x\times S\cup \Gamma_y))/\OO^{\oplus n}(-m\cdot(x\times S\cup \Gamma_y))$$
over $X\times S$, such that the quotient $\OO^{\oplus n}(m\cdot(x\times S\cup \Gamma_y))/{\mathcal J}$
is $S$--flat.

For $m'\geq m$, there is a natural closed embedding $\on{Hilb}^m_n\to\on{Hilb}^{m'}_n$.
It is easy to see that for $G=GL(n)$, $\Gr'_X$ identifies naturally with the inductive limit 
$\underset{\longrightarrow}{\on{Hilb}^m_n}$.

For general $G$, we choose a faithful representation $G\hookrightarrow GL(n)$ 
and show as in Sect. A.5 that $\Gr'_X(G)$ is a closed subfunctor inside $\Gr'_X(GL_n)$.

\end{proof}

\sssec{}

Let $\Fl'_{X\setminus x}$, $\Fl'_x$, $\Gr'_{X\setminus x}$ and $\Gr'_x$ denote the corresponding subschemes
of $\Fl'_X$ and $\Gr'_X$, respectively.

\begin{prop} \label{descr}

There are natural isomorphisms
$$\Fl'_{X\setminus x}\simeq \Gr_{X\setminus x}\times \Fl,\,\,\,
\Gr'_{X\setminus x}\simeq \Gr_{X\setminus x}\times \Gr \text{ and }$$
$$\Fl'_x\simeq \Fl,\,\,\, \Gr'_x\simeq\Gr.$$

\end{prop}

\begin{proof}

The fact that $\Fl'_x\simeq \Fl$ and $\Gr'_x\simeq\Gr$ follows immediately 
from \lemref{bl}. Hence, we must analyze the situation over $X\setminus x$.
We will prove the assertion for $\Fl'_X$, since the proof for $\Gr'_X$
is the same. We will construct canonical morphisms
in both directions between the corresponding functors. 

\smallskip

\noindent $\Rightarrow$  

Let $(y,\F_G,\beta',\epsilon)$ be as above with 
$\Gamma_y\cap x\times S=\emptyset$. First, we define new $G$--bundles $\F^1_G$
and $\F^2_G$ as follows:

$\F^1_G$ (resp., $\F^2_G$) is by definition trivial over 
$X\times S\setminus \Gamma_y$
(resp., $(X\setminus x)\times S$) and is identified with $\F_G$ over
$(X\setminus x)\times S$ (resp., over $X\times S\setminus \Gamma_y$).

Since $(X\setminus x)\times S\cup X\times S\setminus \Gamma_y=X\times S$,
in order to have well-defined $\F^1_G$ and $\F^2_G$ over $X\times S$,
we must define a gluing data over the intersection 
$(X\setminus x)\times S\cap X\times S\setminus \Gamma_y$. However,
the corresponding gluing data for both $\F^1_G$ and $\F^2_G$ are provided
by the isomorphism 
$\beta':\F_G|_{X\times S\setminus (x\times S\cup \Gamma_y)}
\to \F^0_G|_{X\times S\setminus (x\times S\cup \Gamma_y)}$. 

By construction, we have the
trivializations $$\beta^1:\F^1_G|_{X\times S\setminus \Gamma_y}
\to\F^0_G|_{X\times S\setminus \Gamma_y} \text{ and } 
\beta^2:\F^2_G|_{(X\setminus x)\times S}\to \F^0_G|_{(X\setminus x)\times S}.$$

Since $x\times S\in X\times S\setminus \Gamma_y$, the data of $\epsilon$
gives rise to a reduction $\epsilon^2$ of $\F^2_G|_{x\times S}$ to $B$.

Thus, to $(y,\F_G,\beta',\epsilon)$ above we attach the point
$(y,\F^1_G,\beta^1)\times (\F^2_G,\beta^2,\epsilon^2)
\in\Gr_{X\setminus x}\times \Fl$.

\smallskip

\noindent $\Leftarrow$

Let $(y,\F^1_G,\beta^1)\times (\F^2_G,\beta^2,\epsilon^2)$ be an $S$--point 
of $\Gr_{X\setminus x}\times \Fl$. We attach to it a point of 
$\Fl'_X$ as follows:

The $G$--bundle $\F_G$ is by definition identified with $\F^1_G$
over $(X\setminus x)\times S$ and with $\F^2_G$ over 
$X\times S\setminus \Gamma_y$. 

The gluing data for $\F_G$ over $X\times S\setminus \Gamma_y\cap
(X\setminus x)\times S$ is given by the composition:
$$\F^1_G|_{X\times S\setminus \Gamma_y\cap
(X\setminus x)\times S}\overset{\beta^1} \longrightarrow
\F^0_G|_{X\times S\setminus \Gamma_y\cap (X\setminus x)\times S}
\overset{\beta^2}\longleftarrow 
\F^2_G|_{X\times S\setminus \Gamma_y\cap (X\setminus x)\times S}.$$

Thus, we obtain a well-defined $G$--bundle over $X\times S$, which
is trivialized, by construction, over 
$X\times S\setminus \Gamma_y\cap (X\setminus x)\times S$.

Finally, the data of $\epsilon^2$ for $\F^2_G$ defines a data
of $\epsilon$ for $\F_G$, as
$\F_G|_{x\times S}\simeq \F^2_G|_{x\times S}$.

\medskip

Thus, we have constructed maps
$\Fl'_{X\setminus x}\leftrightarrows \Gr_{X\setminus x}\times \Fl$
and it is easy to see that they are inverses of one another.

\end{proof}

\ssec{The functors $\C(\cdot,\cdot)$}

\sssec{}  \label{intrmore}

According to \secref{globalsheaf} and \propref{descr} we can produce a functor
$\on{P}_{\Aut}(\Gr)\times \on{P}(\Fl)\to \on{P}(\Fl'_{X\setminus x})$ by
$$\S,\T\mapsto \S_{X\setminus x}\boxtimes \T$$
and we set
$$\C_{\Fl}(\S,\T):=\Psi_{\Fl'_X}(\S_{X\setminus x}\boxtimes \T).$$

Thus, $\C_{\Fl}(\cdot,\cdot)$ is a functor $\on{P}_{\Aut}(\Gr)\times \on{P}(\Fl)\to \on{P}(\Fl)$. Analogously,
we define the functor $\C_{\Gr}(\cdot,\cdot):\on{P}_{\Aut}(\Gr)\times \on{P}(\Gr)\to \on{P}(\Gr)$ by
setting
$$\C_{\Gr}(\S,\T):=\Psi_{\Gr'_X}(\S_{X\setminus x}\boxtimes \T).$$

\medskip

It is easy to see that $\Fl_X$ is naturally a closed subscheme of $\Fl'_X$:
an $S$--point $(y,\F_G,\beta',\epsilon)$ of $\Fl'_X$ belongs to $\Fl_X$ if and only if the trivialization
$\beta':\F_G|_{X\times S\setminus (x\times S\cup \Gamma_y)}\to\F^0_G|_{X\times S\setminus (x\times S\cup \Gamma_y)}$
extends regularly to $X\setminus \Gamma_y$.

Therefore, we obtain that for $\S\in \on{P}_{\Aut}(\Gr)$,
$$\C_{\Fl}(\S,\delta_{1_{\Fl}})\simeq Z(\S).$$

\medskip

Assertions (a) and (b) of \thmref{main} follow immediately
from the following proposition, whose proof will be given in the next section.

\begin{prop} \label{propone}

Let $\S$ be an object of $\on{P}_{\Aut}(\Gr)$. Then:

\smallskip

\noindent{\rm (a)} For $\T\in \on{P}_{\Iw}(\Fl)$ 
(resp., $\T\in \on{P}_{G(\O)}(\Gr)$)
there is a canonical isomorphism $\C_{\Fl}(\S,\T)\simeq Z(\S)\starfl \T$ 
(resp., $\C_{\Gr}(\S,\T)\simeq\S\stargr\T$).

\smallskip

\noindent{\rm (b)} 
For any $\T\in \on{P}(\Fl)$ (resp., $\T\in \on{P}(\Gr)$)
and $\S\in\on{P}_{G(\O)}(\Gr)$
there is a canonical isomorphism $\C_{\Fl}(\S,\T)\simeq \T\starfl Z(\S)$ 
(resp., $\C_{\Gr}(\S,\T)\simeq\T\stargr\S$).

\end{prop}

Note that the assertion of \propref{propone} for $\Gr$ implies that for $\S\in \on{P}_{G(\O)}(\Gr)$,
$\T\stargr\S$ is perverse for any $\T\in\on{P}(\Gr)$ and that if $\T$ is $G(\O)$--equivariant too, 
then $\S\stargr\T\simeq \T\stargr\S$.

\section{Proofs-I}

\ssec{Some properties of the nearby cycles functor}

\sssec{}

The proof of \thmref{main} will repeatedly use the following well--known result (cf. \cite{SGA}):

Let $\Y$ be a scheme over $X$ and let $\wt{\Y}$ be another scheme with a proper map $g:\wt{\Y}\to \Y$.
Let $g_{X\setminus x}$ (resp., $g_x$) denote the restriction of $g$ to the corresponding subschemes of 
$\wt{\Y}$.

\begin{thm} \label{comdirim}
There is a natural isomorphism of functors $\on{D}^b(\wt{\Y}_{X\setminus x})\to \on{D}^b(\Y_x)$:
$$g_x{}_{!}\circ \Psi_{\wt{\Y}}\simeq \Psi_{\Y}\circ g_{X\setminus x}{}_{!}.$$
\end{thm}

\sssec{}

Let us deduce from \thmref{comdirim} the assertion of of \thmref{main}(d).

\begin{proof}

We will apply \thmref{comdirim} to the map $\pi_X:\Fl_X\to\Gr_X$. We have:

$$\pi_x{}_!(Z(\S))\simeq \pi_{!}(\Psi_{\Fl_X}(\S_{X\setminus x}\boxtimes \delta_{1_{G/B}}))\simeq 
\Psi_{\Gr_X}(\pi_{X\setminus x}{}_{!}(\S_{X\setminus x}\boxtimes \delta_{1_{G/B}})).$$
However, $\pi_{X\setminus x}{}_{!}(\S_{X\setminus x}\boxtimes \delta_{1_{G/B}})\simeq \S_{X\setminus x}$.

Hence, it remains to show that $\Psi_{\Gr_X}(\S_{X\setminus x})\simeq \S$, i.e. that the vanishing cycles functor
$\Phi_{\Gr_X}$ applied to $\S_X$ yields zero. This follows almost immediately from the fact that
$\S_X$ was obtained by the twisted external product construction: 

Let $Y$ be a closed $\Aut$--invariant subscheme of $\Gr$ which contains the support of 
$\S$. Let $\Aut_m$
be a finite dimensional quotient of $\Aut$ such that the action of the latter on $Y$ factors through $\Aut_m$
and let $\X_m$ be the corresponding $\Aut_m$--torsor over $X$. Then $Y_X:=\X_m\overset{\Aut_m}\times Y$ is a closed
subscheme of $\Gr_X$ that contains the support of $\S_X$. Hence, it is enough to calculate $\Phi_{Y_X}(\S_X)$.

\smallskip

The map $\X_m\times Y\to Y_X$ is smooth and it is well-known that the functors
of nearby and vanishing cycles commute with pull-backs under smooth morphisms. Therefore, it is enough to check that
$\Phi_{\X_m\times Y}$ applied to the pull-back of $\S_X$ to $\X_m\times Y$ is $0$. However, the above pull-back 
is a direct product $\Ql{}_{\X_m}\boxtimes \S$. Since the projection $\X_m\to X$ is smooth, this implies the required
vanishing.

\end{proof}

\ssec{Proof of \propref{propone}(a)}

\sssec{}

To prove \propref{propone}(a) we will introduce an auxiliary indscheme $\wt{\Fl'_X}$ over $X$. For a scheme $S$, 
$\Hom(S,\wt{\Fl'_X})$ is the set of 7-tuples $(y,\F_G,\F^1_G,\wt{\beta},\beta^1,\epsilon,\epsilon^1)$, 
where $(y,\F^1_G,\beta^1,\epsilon^1)$ is a point of $\Fl_X$, $\F_G$ is another $G$--bundle, 
$\wt{\beta}$ is an isomorphism $\F_G|_{(X\setminus x)\times S}\to \F^1_G|_{(X\setminus x)\times S}$
and $\epsilon$ is a data of a reduction of $\F_G|_{x\times S}$ to $B$.

By construction, there is a natural projection $p^1_X: \wt{\Fl'_X}\to\Fl_X$ that ``remembers'' only the data of
$(y,\F^1_G,\beta^1,\epsilon^1)$ and a projection $p_X: \wt{\Fl'_X}\to\Fl'_X$ that sends 
$(y,\F_G,\F^1_G,\wt{\beta},\beta^1,\epsilon,\epsilon^1)$ to $(y,\F_G,\beta,\epsilon)$, where
$\beta$ is the composition $\beta^1\circ\wt{\beta}$ defined over $X\times S\setminus (\Gamma_y\cup x\times S)$.

\medskip

By definition, the projection $p^1_X$ makes $\wt{\Fl'_X}$ a fibration over 
$\Fl_X$ with the typical fiber $\Fl$. 
Let us make this assertion more precise. 
Recall that over $\Gr$ we had a ``tower'' of $G(\O)_k$--torsors $\G_k$.
An analogous tower exists globally: 

We introduce a $\Iw_k$--torsor $\G_X{}_k$ over 
$\Fl_X$ that classifies the data of $(y,\F^1_G,\beta^1,\epsilon^1,\gamma_k)$, where
$(y,\F^1_G,\beta^1,\epsilon^1)$ are as in the definition of $\Fl_X$ and
$\gamma_k$ is a data of a trivialization of $\F_G$ on $\D_k\times S$, which is 
compatible with $\epsilon^1$
(i.e. the two reductions to $B$ on $S\simeq\D_0\times S$ coincide).

Let $Y$ be an $\Iw$--invariant closed subscheme of $\Fl$ on which $\Iw$ acts via the quotient $\Iw\to \Iw_k$.
The fibration $Y_X:=\G_X{}_k\overset{I_k}\times Y$
over $\Fl_X$ associated with the $\Iw_k$--torsor $\G_X{}_k$ and the $\Iw_k$--scheme $Y$
is independent of $k$ and 
is naturally a closed ind--subscheme of $\wt{\Fl'_X}$ (the 
latter is an inductive limit of indschemes described in the above way). 

We will denote by $Y_x$ and $Y_{X\setminus x}$ the corresponding subschemes in $Y_X$ (note that $Y_x$
identifies with the corresponding closed sub--indscheme of the convolution diagram $\on{Conv}_{\Fl}$).

\sssec{}

Let $\T$ be as in \propref{propone}(a). Choose $Y$ as above so that 
$\T$ is supported on $Y$. 

As was explained in \secref{intrconstr}, starting from $\S\in \on{P}_{\Aut}(\Gr)$, we can form a perverse sheaf
$\S_{X\setminus x}\boxtimes \delta_{1_{G/B}}$ on $\Fl_{X\setminus x}$ and by taking its twisted external product
with $\T$ we obtain a perverse sheaf 
$(\S_{X\setminus x}\boxtimes \delta_{1_{G/B}})\tboxtimes \T$ on $Y_{X\setminus x}$,
and hence on $\wt{\Fl'_X}$.

Let $p_x$ and $p_{X\setminus x}$ denote the restriction of 
the map $p_X$ to the corresponding subschemes of $\wt{\Fl'_X}$.
The following assertion follows from the definitions:

\begin{lem}

{\em (a)} The map $Y_X\to \wt{\Fl'_X}\overset{p_X}\longrightarrow \Fl'_X$ is proper.

\smallskip

{\em (b)}
The direct image $p_{X\setminus x}{}_{!}((\S_{X\setminus x}\boxtimes \delta_{1_{G/B}})\tboxtimes \T)$ is canonically isomorphic
to the perverse sheaf $\S_{X\setminus x}\boxtimes\T$ on $\Fl'_{X\setminus x}$ constructed in \secref{intrmore}.
\end{lem}

To prove the proposition, let us apply \thmref{comdirim} to the above map
$Y_X\to \wt{\Fl'_X}\overset{p_X}\longrightarrow \Fl'_X$ and the perverse sheaf 
$(\S_{X\setminus x}\boxtimes \delta_{1_{G/B}})\tboxtimes \T$ on $Y_X$. 
The map $p_x:\wt{\Fl'_x}\to \Fl'_x$ identifies with the map
$p:\on{Conv}_{\Fl}\to \Fl$, therefore, it remains to show that 
$\Psi_{Y_X}((\S_{X\setminus x}\boxtimes \delta_{1_{G/B}})\tboxtimes\T)\simeq Z(\S)\tboxtimes \T.$
The argument is similar to the one we used to prove \thmref{main}(d):

The natural projection $\G_X{}_k\times Y\to \G_X{}_k\overset{I_k}\times Y:=Y_X$ is smooth and has
connected fibers. Therefore, it is sufficient to perform the nearby cycles calculation ``upstairs'', i.e.
after the pull-back to $\G_X{}_k\times Y$. However, when we pull-back 
$(\S_{X\setminus x}\boxtimes \delta_{1_{G/B}})\tboxtimes \T$, it decomposes as a direct product 
$(\S_{X\setminus x}\boxtimes \delta_{1_{G/B}})_k\boxtimes \T$, where 
$(\S_{X\setminus x}\boxtimes \delta_{1_{G/B}})_k$ denotes
the pull-back of the perverse sheaf $\S_{X\setminus x}\boxtimes \delta_{1_{G/B}}$ from $\Fl_{X\setminus x}$ to 
$\G_{X\setminus x}{}_k$. 

Let $Z(\S)_k$ denote the pull-back of $Z(\S)$ under the map $\G_x{}_k\to \Fl_x\simeq \Fl$. We have:
$\Psi_{\G_X{}_k}(\S_{X\setminus x}\boxtimes \delta_{1_{G/B}})_k\simeq Z(\S)_k$, since the map $\G_X{}_k\to \Fl_X$
is smooth. Hence,
$$\Psi_{\G_X{}_k\times Y}((\S_{X\setminus x}\boxtimes \delta_{1_{G/B}})_k\boxtimes\T)\simeq Z(\S)_k\boxtimes \T,$$
which is what we had to prove.

\medskip

The proof of \propref{propone}(a) for $\Gr$ is completely similar (and even simpler). 

\ssec{Proof of \propref{propone}(b)}

\sssec{}  \label{beginproofb}

To prove \propref{propone}(b) we will introduce another scheme $\wt{\Fl'_X}$, different from
the one of the previous subsection. (We are going to prove \propref{propone}(b) for $\Fl$, since
the argument for $\Gr$ is the same).

The new $\wt{\Fl'_X}$ classifies 7-tuples 
$(y,\F_G,\F^1_G,\wt{\beta},\beta^1,\epsilon,\epsilon^1)$, 
where $(\F^1_G,\beta^1,\epsilon^1)$ is an $S$--point of $\Fl$ 
(in particular, $\beta^1$ is a trivialization of $\F^1_G$
on $(X\setminus x)\times S$), $\F_G$ is another $G$--bundle, $y$ is another $S$-point of $X$, 
$\wt{\beta}$ is an isomorphism $\F_G|_{X\times S\setminus \Gamma_y}\to 
\F^1_G|_{X\times S\setminus\Gamma_y}$
and $\epsilon$ is a data of a reduction of $\F_G|_{x\times S}$ to $B$.

We have the projections $p^1_X$ and $p_X$ from $\wt{\Fl'_X}$ to $\Fl$ and $\Fl'_X$, respectively:
$p^1_X$ remembers the quadruple $(\F^1_G,\beta^1,\epsilon^1)$ and $p_X$ sends the above
7-tuple to $(y,\F_G,\beta,\epsilon)$, where $\beta$ is the composition $\beta^1\circ\wt{\beta}$ 
defined over $X\times S\setminus (\Gamma_y\cup x\times S)$.

\smallskip

\noindent{\it Remark.}
The essential difference between points (a) and (b) of the proposition is that in the latter
case, $\Fl'_X$ is not strictly speaking a fibration over $\Fl_X$ attached to a group in the sense
of \secref{conventions}. For that reason we have to work harder.

\medskip

Let $\wt{\Fl'_{X\setminus x}}$ 
and $\wt{\Fl'_x}$ denote the corresponding subschemes of $\wt{\Fl'_X}$.
First, observe that 
$\wt{\Fl'_x}$ again identifies canonically with $\Conv_{\Fl}$. Secondly, $\Fl'_{X\setminus x}$
is naturally a closed subscheme in $\wt{\Fl'_{X\setminus x}}$: 

Indeed to an $S$--point 
$(y,\F_G,\beta',\epsilon)$ of $\Fl'_{X\setminus x}$ we attach the data of 
$(y,\F_G,\F^1_G,\wt{\beta},\beta^1,\epsilon,\epsilon^1)$, where 
$(y,\F_G,\epsilon)$ are with no change,
$\F^1_G$ is set to be isomorphic to $\F_G$ over 
$X\times S\setminus \Gamma_y$ and to $\F^0_G$ over
$(X\setminus x)\times S$ (with the gluing data provided by $\beta'$), 
$\epsilon^1$ being induced by $\epsilon$
and $\wt{\beta},\beta^1$ coming by construction.

Thus, for $\S\in \on{P}_{G(\O)}(\Gr)$ and
$\T\in \on{P}(\Fl)$, by taking the direct image of
$\S_{X\setminus x}\boxtimes \T\in \on{P}(\Fl'_{X\setminus x})$, we obtain
a perverse sheaf $(\S_{X\setminus x}\boxtimes \T)^{\sim}$ on $\wt{\Fl'_{X\setminus x}}$.
Its direct image under $p_{X\setminus x}:\wt{\Fl'_{X\setminus x}}\to \Fl'_{X\setminus x}$
is canonically isomorphic to $\S_{X\setminus x}\boxtimes \T$.

By applying \thmref{comdirim} to $(\S_{X\setminus x}\boxtimes \T)^{\sim}$, 
we conclude that it is enough to show that 
$$\Psi_{\wt{\Fl'_X}}((\S_{X\setminus x}\boxtimes \T)^{\sim})\simeq\T\tboxtimes Z(\S).$$

\sssec{}

Let $Y$ be the support of $\T$ in $\Fl$. We can replace our initial $X$ by $X={\mathbb A}^1$ and in the latter
case there exists an \'etale and surjective map $\U\to Y$ such that when we pull-back
the universal $G$--bundle from $X\times\Fl$ to $X\times\U$, 
it becomes trivial. \footnote{Indeed, according to \cite{DS}, we can find $\U$, 
\'etale and surjective over $Y$, such that our $G$-bundle on $X\times \U$ admits a 
reduction to $B$. However, it is well-known that any $B$-bundle on ${\mathbb A}^1\times\U$ is pulled back from
$\U$. Hence, by localizing even more with respect to $\U$, we can make this $B$-bundle trivial.}
Let $\phi'_\U$ be a trivialization;
by further localizing $\U$, we can arrange that the two reductions to 
$B$ of our $G$--bundle on $x\times \U$ 
(one coming from $\phi'_\U$ and the other from the universal property of $\Fl$) coincide.

Let us make a base change $\wt{\Fl'_X}\Rightarrow \wt{\Fl'_X}\underset{\Fl}\times\U$.
Since the projection $\wt{\Fl'_X}\underset{\Fl}\times\U\to \Fl'_X$ is 
\'etale over the support of 
$(\S_{X\setminus x}\boxtimes \T)^{\sim}$, it is enough to perform the nearby cycles 
calculation ``upstairs''.
Namely, let $(\S_{X\setminus x}\boxtimes \T)^{\sim}_\U$ and 
$(\T\tboxtimes Z(\S))_\U$ denote the pull-backs
of $(\S_{X\setminus x}\boxtimes \T)^{\sim}$ and $\T\tboxtimes Z(\S)$ to 
$\wt{\Fl'_{X\setminus x}}\underset{\Fl}\times\U$ and 
$\wt{\Fl'_x}\underset{\Fl}\times\U$, respectively.

We must show that
$$\Psi_{\wt{\Fl'_X}\underset{\Fl}\times\U}
((\S_{X\setminus x}\boxtimes \T)^{\sim}_\U)\simeq (\T\tboxtimes Z(\S))_\U$$
and that this isomorphism is independent of the choice of the trivialization $\phi'_\U$.

\medskip

Notice now, that the choice of $\phi'_\U$ defines an
identification $\wt{\Fl'_X}\underset{\Fl}\times\U\simeq \Fl_X\times \U$.
Indeed, since the univesral bundle $\F'_G$ over $X\times \U$ is trivial, the 
data of $\wt{\beta}$ is
equivalent to the trivialization of $\F_G$ off $\Gamma_y$.

However, when we restrict the universal bundle to $(X\setminus x)\times\U$, 
it has two different trivializations!
One comes from $\phi'_\U$ and the other from the fact that the universal bundle on 
$X\times\Fl$ is by definition 
trivialized over $(X\setminus x)\times\U$. 
These trivializations differ by a map $\phi_\U:(X\setminus x)\times\U\to G$. 

\medskip

Correspondingly, we have two different closed embeddings of $\Gr_{X\setminus x}\times \U$
into $\wt{\Fl'_{X\setminus x}}\underset{\Fl}\times\U$:

\noindent Embedding (1) is the composition of the above identification 
$\wt{\Fl'_X}\underset{\Fl}\times\U\simeq \Fl_X\times \U$ and the embedding
$\Gr_{X\setminus x}=\Gr_{X\setminus x}\times 1_{G/B}\hookrightarrow
\Fl_{X\setminus x}$ of \propref{descrFLX}.

\noindent Embedding (2) comes from the embedding 
$\Fl'_{X\setminus x}\to \wt{\Fl'_{X\setminus x}}$
described in \secref{beginproofb} and the isomorphism 
$\Fl'_{X\setminus x}\simeq \Gr_{X\setminus x}\times\Fl$
of \propref{descr}.

It is easy to see that these two embeddings differ by the 
automorphism of $\Gr_{X\setminus x}\times \U$ induced by $\phi_\U$ 
as in \propref{expequiv}. 

\medskip

By construction, $(\S_{X\setminus x}\boxtimes \T)^{\sim}_\U$ is isomorphic to the direct image
under the above Embedding (2) of $\S_{X\setminus x}\boxtimes \T_\U$ (here $\T_\U$
is the pull-back of $\T$ under $\U\to\Fl$). Hence, by 
\propref{expequiv}, it is isomorphic also to the direct image of the same 
$\S_{X\setminus x}\boxtimes\T_\U$
under Embedding (1). 

Hence, on the one hand,
$$\Psi_{\wt{\Fl'_X}\underset{\Fl}\times\U}
((\S_{X\setminus x}\boxtimes \T)^{\sim}_\U)\simeq 
\Psi_{\Fl_X\times \U}((\S_{X\setminus x}\boxtimes\delta_{1_{G/B}})
\boxtimes \T_\U)\simeq Z(\S)\boxtimes\T_\U.$$

But on the one hand, under the identification
$\wt{\Fl'_x}\underset{\Fl}\times\U\simeq \Fl\times\U$, 
the complex $(\T\tboxtimes Z(\S))_\U$ goes over to the same $Z(\S)\boxtimes\T_\U$.

This proves the existence of the required isomorphism. Let us now analyze what happens when 
we modify $\phi'_\U$ by a map $\phi''_\U:X\times \U\to G$ 
($\phi''_\U$ must send $x\times \U$ to $B\subset G$).
The effect would be the automorphism of $\Fl_X\times \U$ 
induced by $\phi''_\U$, as in \propref{expequiv}. This does not change
the identification $\Psi_{\wt{\Fl'_X}\underset{\Fl}\times\U}
((\S_{X\setminus x}\boxtimes \T)^{\sim}_\U)\simeq (\T\tboxtimes Z(\S))_\U$, by the 
definition of the $\Iw$-equivariant structure on $Z(\S)$ (cf. proof of \propref{Iwequiv}).

\section{Proofs-II}

\ssec{The monodromy action}

\sssec{}

Recall the situation of \secref{nearby}. Let $\Gamma$ (resp., $\Gamma^g$) denote the full 
(resp., geometric) Galois group that corresponds to the pair $x\in X$. 
In other words, $\Gamma$ (resp., $\Gamma^g$) is the Galois group of the field of 
fractions of the
henselization (resp., strict henselization) of the local ring $\OO_x$. 
As was mentioned in \secref{nearby}, we are fixing a splitting  
$\Gamma\simeq \Gamma^g\rtimes\widehat{\on{Gal}(\Fqb/\Fq)}$.
There is a canonical homomorphism $\Gamma^g\overset{t_\ell}\to {\mathbb Z}_\ell(1)$
(we are taking into account the action of $\on{Gal}(\Fqb/\Fq)$ on 
$\Gamma^g$ and ${\mathbb Z}_\ell(1)$).

Let $\rho:\Gamma\to \on{Aut}(V)$ be a (continuous) representation.
Following \cite{Gr}, there exists a canonical nilpotent endomorphism
$M_V:V\to V(-1)$ and a subgroup $\Gamma'\subset\Gamma^g$ of finite index,
such that for any $\gamma\in \Gamma'$
$$\rho(\gamma)=\on{exp}(t_\ell(\gamma)\cdot M_V):V\to V.$$

Recall that a representation $(\rho,V)$ is called unipotent if $\Gamma'=\Gamma^g$.
In this case, there exists a $\Gamma$-stable filtration on $V$ such that
the action of $\Gamma^g$ on the successive quotients is trivial.

In general, any representation $\rho:\Gamma\to \on{Aut}(V)$ can be decomposed as a direct 
sum $V=V^{un}\oplus V^{non-un}$,
where $V^{un}$ is unipotent and $V^{non-un}$ is purely non--unipotent (i.e. every irreducible 
subquotient of $V^{non-un}$ is non-trivial as a $\Gamma^g$--representation).

\medskip

Now, the basic property of the functor $\Psi_\Y$ is that it carries the action of $\Gamma$. 
In particular, for $\S\in\on{P}(\Y_{X\setminus x})$, 
we have a nilpotent endomorphism $M_\S:\Psi_\Y(\S)\to \Psi_\Y(\S)(-1)$ and a 
decomposition $\Psi_\Y(\S)\simeq \Psi^{un}_\Y(\S)\oplus \Psi^{non-un}_\Y(\S)$.

\begin{lem} \label{commore}
The isomorphism of functors given by \thmref{comdirim} respects the action of $\Gamma$.
\end{lem}

\sssec{}

Let us apply the above discussion to the situation $\Y=\Fl_X$. We obtain for every 
$\S\in\on{P}_{G(\O)}(\Gr)$ an endomorphism $M_\S:Z(\S)\to Z(\S)(-1)$ and 
a canonical decomposition $Z(\S)=Z(\S)^{un}\oplus Z(\S)^{non-un}$.

The following is not essential for our purposes, but is important as an observation:

\begin{prop} \label{unipotent}
$Z(\S)^{non-un}=0$.
\end{prop}

\begin{proof}

First, by the construction of the functor $\C_{\Fl}(\cdot,\cdot)$, we have a 
$\Gamma$--action on 
$\C_{\Fl}(\S,\T)$ for every $\T\in\on{P}(\Fl)$. 
Moreover, \lemref{commore} implies that the isomorphisms 
$\C_{\Fl}(\S,\T)\simeq \T\starfl Z(\S)$ and
$\C_{\Fl}(\S,\T)\simeq Z(\S)\starfl\T$ (for $\T\in\on{P}_{\Iw}(\Fl)$) 
are compatible with this $\Gamma$--action. 

Therefore, for $\T\in\on{P}_{\Iw}(\Fl)$ we have the isomorphisms
$$Z(\S)^{un}\starfl\T\simeq \T\starfl Z(\S)^{un} \text{ and } 
Z(\S)^{non-un}\starfl\T\simeq \T\starfl Z(\S)^{non-un}.$$

In addition, \lemref{commore} applied to the map $\pi_X:\Fl_X\to\Gr_X$ implies that
$\pi_{!}(Z(\S)^{non-un})\simeq\Psi^{non-un}_{\Gr_X}(S_{X\setminus x})=0$.

\medskip

Let us study the function corresponding to $Z(\S)^{non-un}$ on $\Fl(\Fq)$. (By enlarging
the finite field we may assume that $G$ is split.)

The fact
that $Z(\S)^{non-un}\starfl\T$ and $\T\starfl Z(\S)^{non-un}$ are isomorphic as perverse 
sheaves defined over $\Fq$
for any $\T\in\on{P}_{\Iw}(\Fl)$, implies that the corresponding element of 
$\HH_{\Iw}$ is central.
At the same time, it vanishes under the map $\pi:\HH_{\Iw}\to\HH_{sph}$. Hence, this function 
is zero, as $\pi$ induces an isomorphism $Z(\HH_{\Iw})\to\HH_{sph}$, by Bernstein's theorem.

Now, the same fact is true not only for $\Fq$, but also for all
finite field extensions $\Fq\subset{\mathbb F}_{q'}$, which implies that $Z(\S)^{non-un}=0$.

\end{proof}

\noindent{\it Remark.}
The above proof of \propref{unipotent} uses the ``faisceaux-fonctions''
correspondence and Bernstein's theorem. In fact, it is not difficult to give
a purely geometric proof, which we will do elsewhere. 

\ssec{Proof of \thmref{main}(c)}

\sssec{}

The proof will rely on the following general property of the nearby cycles 
functor, proved in \cite{BB}:

Now let $\Y^1$ and $\Y^2$ be two schemes mapping to $X$ and let $\S^1$ and $\S^2$ be objects in 
$\on{D}^b(\Y^1_{X\setminus x})$
and $\on{D}^b(\Y^2_{X\setminus x})$, respectively. We will denote by 
$\S^1\underset{X}\boxtimes\S^2$ the $*$-restriction of $\S^1\boxtimes\S^2$ to 
$\Y^1_{X\setminus x}\underset{X\setminus x}\times \Y^2_{X\setminus x}$, cohomologically shifted by $1$ to the right.

\begin{thm} \label{bb}
There is a canonical isomorphism in $\on{D}^b(\Y^1_x\times\Y^2_x)$:
$$\Psi_{\Y^1\underset{X}\times \Y^2}(\S^1\underset{X}\boxtimes\S^2)\simeq 
\Psi_{\Y^1}(\S^1)\boxtimes \Psi_{\Y^2}(\S^2).$$
Moreover, this isomorphism is compatible with the $\Gamma$--action.
\end{thm}

\sssec{}

To prove \thmref{main}(c), we introduce the schemes $\Conv_X$ and $\Conv'_X$ over $X$. 
The scheme $\Conv_X$ is by definition $\X\underset{\Aut}\times\Conv_{\Gr}$. 
In other words, $\Hom(S,\Conv_X)$
is the set of quintuples $(y,\F_G,\wt{\beta},\F^1_G,\beta^1)$, where 
$(y,\F^1_G,\beta^1)$ is a point of $\Gr_X$, $\F_G$ is another $G$--bundle over $X\times S$ and
$\wt{\beta}$ is an isomorphism $\F_G|_{X\times S\setminus\Gamma_y}\simeq 
\F^1_G|_{X\times S\setminus\Gamma_y}$.

For $\Conv'_X$, we put $\Hom(S,\Conv'_X)$ to be the set 
of 7-tuples $(y,\F_G,\wt{\beta},\epsilon,\F^1_G,\beta^1,\epsilon^1)$, where 
$(y,\F_G,\wt{\beta},\F^1_G,\beta^1)$ are as above and $\epsilon$ (resp., $\epsilon^1$)
is a reduction to $B$ of $\F_G|_{x\times S}$ (resp., of $\F^1_G|_{x\times S}$).

Let $\Conv_{X\setminus x}$, $\Conv_x$, $\Conv'_{X\setminus x}$ and $\Conv'_x$ denote the corresponding subschemes of 
$\Conv_X$ and $\Conv'_X$, respectively. Let also $p_X$ and $p^1_X$ denote the standard projections from $\Conv'_X$ to $\Fl_X$. 

\begin{lem} \label{descrc}
We have natural identifications
$$\Conv'_{X\setminus x}\simeq \Conv_{X\setminus x}\times G/B\times G/B \text{ and }
\Conv'_x\simeq\Conv_{\Fl}.$$
\end{lem}

Thus, starting from two objects $\S^1$ and $\S^2$ of $\on{P}_{G(\O)}(\Gr)$, we can construct a perverse sheaf
$\S^1_{X\setminus x}\tboxtimes \S^2_{X\setminus x}\boxtimes \delta_{1_{G/B}}\boxtimes \delta_{1_{G/B}}$
on $\Conv'_{X\setminus x}$.
By applying \thmref{comdirim} to $p_X:\Conv'_X\to\Fl_X$, we obtain that in order to prove
\thmref{main}(c) and \thmref{monodromy}, we must verify the following:

There exists a $\Gamma$--equivariant isomorphism
$$\Psi_{\Conv'_X}(\S^1_{X\setminus x}\tboxtimes \S^2_{X\setminus x}\boxtimes \delta_{1_{G/B}}\boxtimes \delta_{1_{G/B}})\simeq
Z(\S^1)\tboxtimes Z(\S^2) \in\on{P}(\Conv_{\Fl}).$$ 

This proof of this statement is a variation of the argument presented in the proof of \propref{propone}(b). We allow
ourselves to be more sketchy:

\sssec{}

Let $Y$ be the support of 
$\S^1_{X\setminus x}\boxtimes\delta_{1_{G/B}}$ in $\Fl_X$. As in the proof of 
\propref{propone}(b), we can assume that
there exists a surjective \'etale map $\U\to Y$ with the property that
the pull-back of the universal $G$--bundle from $X\times \Fl_X$ to $X\times \U$ is trivial. 
Let us fix a trivialization compatible with the existing $B$--structure on $x\times \U$.
Let $\U_{X\setminus x}$ and $\U_x$ denote the corresponding subschemes of $\U$.

As before, it sufficient to carry out the nearby cycles calculation on
$\U\underset{\Fl_X}\times\Conv'_X$, where the fiber product is defined using the projection 
$p^1_X: \Conv'_X\to \Fl_X$.

\smallskip

Let $\E(\S^1,\S^2)$ denote the pull-back of the perverse sheaf
$\S^1_{X\setminus x}\tboxtimes \S^2_{X\setminus x}\boxtimes 
\delta_{1_{G/B}}\boxtimes \delta_{1_{G/B}}$
from $\Conv'_{X\setminus x}$ to 
$\U_{X\setminus x}\underset{\Fl_{X\setminus x}}\times\Conv'_{X\setminus x}$.
Let $\E(\S^1)_{X\setminus x}$ (resp., $\E(\S^1)_x$) denote the pull-back of 
$\S^1_{X\setminus x}\boxtimes\delta_{1_{G/B}}$ 
(resp., of $Z(\S_1)$) to $\U_{X\setminus x}$ (resp., to $\U_x$).

As in the proof of \propref{propone}(b), 
the trivialization of the pulled-back universal $G$--bundle
on $X\times \U$ defines an isomorphism
$\U\underset{\Fl_X}\times\Conv'_X\simeq \U\underset{X}\times \Fl_X$. 
Moreover, under this isomorphism
the perverse sheaf $\E(\S^1,\S^2)$ becomes identified with
$\E(\S^1)_{X\setminus x}\underset{X}\boxtimes (\S^2_{X\setminus x}\boxtimes\delta_{1_{G/B}})$.

Therefore, on the one hand, using \thmref{bb}, we obtain 
$$\Psi_{\U\underset{\Fl_X}\times\Conv'_X}(\E(\S^1,\S^2))\simeq
\Psi_{\U}(\E(\S^1)_{X\setminus x})
\boxtimes \Psi_{\Fl_X}(\S^2_{X\setminus x}\boxtimes\delta_{1_{G/B}})\simeq
\E(\S^1)_x\boxtimes Z(\S_2).$$

On the other hand, the pull-back of 
$Z(\S_1)\tboxtimes Z(\S_2)$ under $\U_x\times \Fl_x\to \Conv_{\Fl}$ 
identifies also with $\E(\S^1)_x\boxtimes Z(\S_2)$, which is what we had to prove.

\section{Appendix} \label{App}

\noindent{\bf A.1.}
Let $F$, $F'$ be two countravariant functors $\on{Schemes}\to\on{Sets}$
and let $F'\to F$ be a morphism. We say that $F'$ is a closed subfunctor of $F$ if for 
any scheme $S$ and any $f_S\in F(S)$ 
the Cartesian product functor $\underline{S}\underset{F}\times F'$
is representable by a closed subscheme of $S$.

Let $Y_1\hookrightarrow Y_2 \hookrightarrow...\hookrightarrow Y_n\hookrightarrow...$
be a directed system of schemes, where all maps $Y_i\hookrightarrow Y_{i+1}$ are closed
embeddings.

We define a functor $\underset{\longrightarrow}{Y_i}$ on the category of schemes by setting
$\on{Hom}(S,\underset{\longrightarrow}{Y_i}):=\underset{\longrightarrow}{\on{Hom}}(S,Y_i)$
for $S$ quasi-compact, and by extending it to all schemes by requiring that it is a sheaf in
Zarisky topology.

A (strict) indscheme is, by definition, a functor $F$ which is isomorphic to some 
$\underset{\longrightarrow}{Y_i}$ as above. We say that an indscheme is of ind--finite type
if the above family of $Y_i$'s can be chosen in such a way that all of them are of finite type.

In what follows we will work with indschemes of ind--finite type only. 
The basic objects of this paper, that
is $\Gr_X$, $\Fl$, $\Gr_X$, etc., all have this property. 
The only indscheme not of ind--finite type
that appears in this paper is $G(\K)$, but it has been used only as a functor. 
Thus, unless specified otherwise, 
by an indscheme we will mean an indscheme of ind--finite type.

\begin{lem} \label{indcorr}
Let $F\simeq \underset{\longrightarrow}{Y_i}$ be an indscheme. Then:

\smallskip

\noindent{\em (a)}
If $Z$ is a scheme of finite type and $\underline{Z}\to F$ is a closed subfunctor, then
there exists an index $i$ such that $Z$ is a closed subscheme of $Y_i$. 
(In this case we will say that
$Z$ is a closed subscheme of $F$.)

\smallskip

\noindent{\em (b)}
If $F'$ is a closed subfunctor of $F$, then $F'$ is an indscheme too.

\smallskip

\noindent{\em (c)}
If $F$ is isomorphic to $\underset{\longrightarrow}{Y'_i}$ 
for a different family of schemes $Y'_i$,
then for every $i_1$ there exists an $i_2$ 
such that $Y'_{i_1}$ is a closed subscheme of $Y_{i_2}$
and vice versa.
\end{lem}

The proof is a tautology. 

\medskip

\noindent{\bf A.2.} 
For an indscheme $F=\underset{\longrightarrow}{Y_i}$ we define the category of perverse
sheaves on it as $\on{P}(F):=\underset{\longrightarrow}{\on{P}}(Y_i)$, where the functors 
$\on{P}(Y_i)\to \on{P}(Y_{i+1})$
are, of course, the direct image functors. This is again an abelian category, since 
the ``direct image under a closed embedding'' functor is exact. Similarly, one can 
define the derived category
$\on{D}^b(F):=\underset{\longrightarrow}{\on{D}^b}(Y_i)$, which is a 
triangulated category due to the 
exactness property mentioned above. Actually, we do no need derived 
categories in this paper and we discuss them only
for the sake of completeness.

\smallskip

Point (3) of \lemref{indcorr} 
implies that these definitions do not depend on the choice of a presentation
of $F$ as $\underset{\longrightarrow}{Y_i}$, 
i.e. $\on{P}(F)$ and $\on{D}^b(F)$ are intrinsically attached to $F$.

We emphasize again, that a perverse sheaf on an indscheme 
is by definition supported on a closed subscheme
of finite type. This means that this notion is essentially ``finite-dimentional''.

\medskip

\noindent{\bf A.3.} 
Let $H$ be a group--scheme, which is a projective limit of linear algebraic groups:
$H=\underset{\longleftarrow}{H_k}$. The basic examples are $H=G(\O)$, $H=\Aut$.

\smallskip

Let $F$ be an indscheme and let $\underline{H}$ act on $F$ (in the sense of functors). 
We say that this action is ``nice'' if the following holds:
every closed subscheme $Z$ of $F$ 
is contained in a larger closed subscheme $Z'$ with the property that $Z'$
is $H$--stable and the action of $H$ on $Z'$ factors through some $H_k$.

\smallskip

Let $F=\underset{\longrightarrow}{Y_i}$ be an indscheme and $H$ a group--scheme of 
the above type. An $H$--torsor 
$\H$ over $F$ is by definition a compatible system of $H_k$--torsors $\H_{i,k}$
over the $Y_i$'s. Again, by point (3) of \lemref{indcorr}, 
this notion does not depend on the presentation
of $F$ as an inductive limit. 
For an $H$--torsor $\H$ we could consider its total space, which will be
an indscheme not of ind--finite type (unless $H$ is finite dimensional), but this will not be of use for us.

Consider the following situation: let $\H$ be an $H$--torsor over $F$ and let $F'$ be another ind-scheme with a
``nice'' action of $H$. We claim, that we can form the ``associated bundle'' $\H\overset{H}\times F'$ over $F$,
which will be again an indscheme.

Indeed, we can represent $F'$ as an inductive limit $\underset{\longrightarrow}{Y'_i}$ with each $Y'_k$ being
an $H_k$--scheme. We set $\H\overset{H}\times F'$ to be the inductive limit of
$\H_{i,i}\overset{H_i}\times Y'_i$.

\medskip

\noindent{\bf A.4.}
Now we will introduce the category of $H$--equivariant perverse sheaves on an indscheme with 
a ``nice'' action of $H$. First, let $Y$ be a scheme of finite type 
and $H'$ an algebraic group acting on it. Then the notion
of an $H'$--equivariant perverse sheaf on $Y$ is well-known. 
If $H''\to H'$ be a surjection, then the categories
$\on{P}_{H'}(Y)$ and $\on{P}_{H''}(Y)$ are naturally equivalent.

Therefore, if $F=\underset{\longrightarrow}{Y_i}$ with each $Y_k$ being 
$H$--stable and acted on via $H\to H_k$,
we can define $\on{P}_H(F)$ as $\underset{\longrightarrow}{\on{P}_{H_i}}(Y_i)$.

In particular, let $\T$ be a perverse sheaf on $F$, $\H$ be an $H$-torsor over $F$,
$F'$ be an indscheme with a ``nice'' $H$-action and, finally, let $\S$ be
an $H$-equivariant perverse sheaf on $F'$. Then the construction of the 
twisted external product (cf. \secref{conventions}) makes sense and it produces
a perverse sheaf $\T\tboxtimes\S$ on $\H\overset{H}\times F'$.

To introduce the derived categories, we need to make an assumption that $H$ contains a subgroup of 
finite codimension which is pro-unipotent, i.e. that for large enough $k$, $H^k:=\on{ker}(H\to H_k)$ 
is a projective limit of unipotent groups. (This assumption is valid in our examples.)

If $Y$ be a scheme of finite type and $H'$ an algebraic group acting on it, the derived category
$\on{D}_{H'}^b(Y)$ has been introduced by Bernstein and Lunts in \cite{BeLu}. Their definition has the following
property: if $H''\to H'$ is a surjection with a unipotent kernel, then the categories 
$\on{D}_{H''}^b(Y)$ and $\on{D}_{H''}^b(Y)$ 
are equivalent. This enables us to introduce $\on{D}^b_H(F)$: by shifting the indices, we may assume that all the
appearing $H^k$'s are already pro-unipotent, and we set 
$\on{D}^b_H(F)=\underset{\longrightarrow}{\on{D}^b_{H_i}}(Y_i)$.

\medskip

\noindent{\bf A.5.}
Recall our definition of $\Gr$ given in \secref{grfl}. 
Here, for the sake of completeness, we will
prove that $\Gr$ is indeed representable by an indscheme. 
Moreover, from the proof it will follow that
any closed subscheme of $\Gr$ is proper over $\Fq$.  We proceed in two steps:

\begin{proof} $   $

\noindent{\bf Step 1.}
Let us first prove the assertion when $G=GL(n)$. In this case for an affine $S$,
$\on{Hom}(S,\Gr)$ consists of pairs $(\V,\beta^V)$, 
where $\V$ is a projective rank $n$ module over
$\OO_S[[t]]$ and $\beta^V$ is an isomorphism 
$\V\underset{\Fq[[t]]}\otimes \Fq((t))\simeq V\otimes\OO_S((t))$, where
$V$ is the standard $n$-dimensional vector space.

For an integer $m$ consider the 
$(2m+1)\cdot n$--dimensional vector space $t^{-m}V[[t]]/t^{m+1}V[[t]]$ and 
consider the functor that associates to a scheme $S$ the set 
of all $S$--flat and $t$--stable submodules
of $\OO_S\otimes t^{-m}V[[t]]/t^{m+1}V[[t]]$. 
This functor is representable by a closed subscheme,
call it $\Gr^m$, of the Grassmannian of $t^{-m}V[[t]]/t^{m+1}V[[t]]$. 
In particular, $\Gr^m$ is proper.

\smallskip

We have a natural closed embedding $\Gr^i\subset \Gr^{i+1}$, 
since $t^{-i}V[[t]]/t^{i+1}V[[t]]$
is canonically a $t$--invariant subquotient of $t^{-i-1}V[[t]]/t^{i+2}V[[t]]$. 
Hence, we obtain a directed
family $\Gr^1\hookrightarrow...\hookrightarrow 
\Gr^i\hookrightarrow \Gr^{i+1}\hookrightarrow...$,
and it is clear that the functor $\Gr$ is isomorphic to $\underset{\longrightarrow}{\Gr^i}$.

\smallskip

\noindent{\bf Step 2.}
To treat the case of an arbitrary $G$, it suffices to show that if $G_1\to G_2$ is an embedding of
reductive groups, then the natural map $\Gr(G_1)\to \Gr(G_2)$ realizes $\Gr(G_1)$ as a closed
subfunctor of $\Gr(G_2)$. For the proof, it will be more convenient to use the realization of the affine Grassmannian via
a curve $X$, instead of the formal disc (cf. \lemref{bl}).
 
Consider the quotient $G_2/G_1$. 
It is a basic fact that under the above circumstances,
$G_2/G_1$ is an affine variety; let $1_{G_2/G_1}\in G_2/G_1$ denote
the point corresponding to $1\in G_2$.

Thus, let $S$ be a scheme and let $\F_{G_2}$ be a $G_2$-bundle on $X\times S$, trivialized by means of $\beta$
over $(X\setminus x)\times S$. In particular, we obtain a map 
$$(X\setminus x)\times S\overset{\beta(1_{G_2/G_1})}\longrightarrow \F_{G_2}\overset{G_2} \times G_2/G_1.$$

Since $\F_{G_2}\overset{G_2} \times G_2/G_1$ is affine over $X\times S$,
there exists a closed subscheme $S'$ of $S$, such that for any $\phi:S''\to S$ such that
$(X\setminus x)\times S''\to \F_{G_2}\overset{G_2} \times G_2/G_1$
extends to a map $X\times S''\to \F_{G_2}\overset{G_2} \times G_2/G_1$, $\phi$ factors 
as $S''\to S'\to S$.

It is easy to see that this $S'$ represents the fiber product
$S\underset{\Gr(G_2)}\times \Gr(G_1)$.

\end{proof}

The above argument proves the representability of $\Gr(G)$ for a reductive group $G$. However, $\Gr(G)$
is representable for an arbitrary linear algebraic group: 

One shows that if 
$$1\to {\mathbb G}_a\to G_1\to G_2\to 1$$
is a short exact sequence of algebraic groups and $\Gr(G_2)$ is representable, then $\Gr(G_1)$ is represenatble too.
For example, when $G_2=\{1\}$, $\Gr({\mathbb G}_a)$ is isomorphic to the direct limit of linear spaces 
$\underset{\longrightarrow} {t^{-i}\cdot \Fq[[t]]/\Fq[[t]]}$ and the general case is not much different.

However, the above example shows that when $G$ is not reductive, $\Gr$ is not a limit of proper (i.e. compact)
schemes. In general, it is easy to see that for an embedding of algebraic groups $G_1\to G_2$, the map $\Gr(G_1)\to \Gr(G_2)$ 
is a locally-closed embedding.

\medskip

\noindent{\bf A.6.} Finally, let us prove \lemref{functconv}. 
This should be well-known and we include
the proof just in order to demonstrate ``how things work''.

On the one hand, we have the indscheme 
$\G\overset{G(\O)}\times\Gr$ (cf. Sect. A.3 above) and on the other hand 
the functor that associates to a scheme $S$ the set of quadruples 
$(\F_G,\F^1_G,\wt{\beta},\beta^1)$, as in
the formulation of the lemma. Both are sheaves on the category
of schemes with the \'etale topology. 

Thus, let $(\F^1_G,\beta^1)$ be an $S$--point of $\Gr$ and we must show that the 
additional data of $(\F_G,\wt{\beta})$ is equivalent to choosing an $S$--point of 
$\G\overset{G(\O)}\times\Gr$ that projects to our point of $\Gr$. 

By making an \'etale localization, we can assume that $\F^1_G$ can be trivialized and 
let us choose such a
trivialization. Then, on the one hand, the data of $(\F_G,\wt{\beta})$ 
becomes equivalent to a data
of just another $S$--point of $\Gr$.

Similarly, once $\F^1_G$ is trivialzed, the induced $G(\O)$--torsor over 
$\G\underset{\Gr}\times S$
becomes trivialized too (which means that all the $G(\O)_k$--torsors become 
trivialized in a compatible way).
Hence $S\underset{\Gr}\times (\G\overset{G(\O)}\times\Gr)$ splits as $S\times \Gr$. 

This proves the required assertion, once we check that the change of a trivialization of
$\F^1_G$ has the same effect on both sides, which is obvious.

\medskip

\noindent{\bf A.7.}
This subsection is not logically related to anything else in the this appendix. We will reprove
the Kazhdan-Lusztig parity vanishing assertion, which was used in the proof of \propref{GOeq}.
The argument presented below is presumably well-known:

We will prove a more general fact, namely, the parity vanishing for stalks of IC sheaves of
$\Iw$-orbits on $\Fl$. For an element $\on{w}$ in the affine Weyl group, let $\Fl_{\on{w}}$
denote the corresponding $\Iw$-orbit in $\Fl$ and $\ol{\Fl}_{\on{w}}$ its closure.
Let $\on{y}$ be another element with $\Fl_{\on{y}}\subset \ol{\Fl}_{\on{w}}$. We must show
that the complex $\IC_{\ol{\Fl}_{\on{w}}}|_{\Fl_{\on{y}}}$ has cohomology (in either usual or perverse sense)
in degrees of constant parity.

In the beginning of the paper we introduced the scheme $\Conv_{\Fl}$, which was a ``two-fold''
convolution of $\Fl$ with itself. Along the same lines one introduces the $n$-fold convolution
$\Conv^n_{\Fl}$. We will denote by $p_n$ the natural projection $\Conv^n_{\Fl}\to \Fl$,
which generalizes the projection $p$, when $n=2$. If $\on{w}_1,...,\on{w}_n$ are the $n$ elements 
of the affine Weyl group, we will denote the corresponding closed subscheme of $\Conv^n_{\Fl}$ by 
$\Conv^{\on{w}_1,...,\on{w}_n}_{\Fl}$.

Now, let $\on{w}=s_1\cdot ...\cdot s_n$ be a reduced decomposition of $\on{w}$ as a product of simple reflections.
Then the restriction of $p_n$ to $\Conv^{s_1,...,s_n}_{\Fl}$ is a proper dominant map onto $\ol{\Fl}_{\on{w}}$.
By the decomposition theorem, $\IC_{\ol{\Fl}_{\on{w}}}$ is a direct summand of 
$p_n{}_!(\IC_{\Conv^{s_1,...,s_n}_{\Fl}})$. Therefore, it suffices to prove the parity vanishing for the fibers
of $p_n{}_!(\IC_{\Conv^{s_1,...,s_n}_{\Fl}})$.

However, since each $\ol{\Fl}_{s_i}\simeq \PP^1$, the variety
$\Conv^{s_1,...,s_n}_{\Fl}$ is non-singular (it is commonly referred to as the Bott-Samelson
resolution of $\ol{\Fl}_{\on{w}}$), hence, $\IC_{\Conv^{s_1,...,s_n}_{\Fl}}$ is the constant sheaf,
up to a cohomological shift. By base change, it therefore suffices to show that the fiber of
$\Conv^{s_1,...,s_n}_{\Fl}$ over every given point in $\ol{\Fl}_{\on{y}}$ has cohomology in degrees
of constant parity. 

Now, it is known (and easily proven by induction) that each such fiber, call it $Y$, can be represented
as is a union on locally closed subvarieties $Y_i$, each of which is isomorphic to a tower of
affine spaces. Therefore, the Cousin spectral sequence implies that $H_c(Y)$ lives only in even degrees.

\end{document}